\PassOptionsToPackage{numbers,sort&compress}{natbib}

\documentclass{article}
\usepackage{amsmath, amssymb, amsthm, graphicx}
\usepackage{hyperref}
\usepackage{natbib}
\usepackage{float}
\graphicspath{{figs/}{./}}
\usepackage{booktabs}
\usepackage{microtype}

\newtheorem{definition}{Definition}[section]
\newtheorem{theorem}{Theorem}[section]
\newtheorem{lemma}{Lemma}[section]
\newtheorem{remark}{Remark}[section]
\newtheorem{corollary}{Corollary}[section]
\newtheorem{proposition}{Proposition}[section]
\newcommand{\keywords}[1]{\textbf{Keywords:} #1}
\usepackage[inner=2.5cm,outer=2.5cm,top=2.5cm,bottom=2.5cm]{geometry}

\usepackage[skip=4pt]{caption}

\newcommand{\betaeta}{\beta_\eta}
\newcommand{\peps}{p_\epsilon}
\newcommand{\Rz}{R_0}
\newcommand{\Rs}{R_0^{s}}
\newcommand{\Lext}{\Lambda_{\mathrm{ext}}}
\newcommand{\Lsharp}{\Lambda_{\mathrm{ext}}^{\dagger}}
\newcommand{\avg}[1]{\big\langle #1 \big\rangle_t}

\title{Stochastic Dynamics of Hepatitis B Virus Infection: Extinction, Persistence, and Ergodicity}

\author{
 Abdallah Alsammani\\
  \textit{Department of Mathematical Sciences, Delaware State University, Dover, DE, USA} \\
  aalsammani@desu.edu
}

\date{}
\begin{document}
\maketitle

\begin{abstract}
\noindent We develop and analyze a three-compartment stochastic model of within-host hepatitis B virus infection with multiplicative It\^o noise. We establish global well-posedness, strict positivity, and uniform moment bounds, and derive an explicit sufficient criterion for almost sure exponential clearance. Standard extinction arguments based on $\ln(y+z)$ are conservative and fail to recover the deterministic spectral threshold. To overcome this limitation, we optimize over the weighted family $\ln(\alpha y+z)$ and control hepatocyte fluctuations through exact occupation identities for $x+y$, avoiding both equilibrium substitution and an invalid pathwise comparison. The resulting general extinction bound admits a reduced optimized form under an explicit admissibility condition. This reduced bound lies no more than $\tfrac12\max\{\sigma_2^2,\sigma_3^2\}$ below the spectral abscissa $\lambda_{\max}(J)$ of the deterministic infected subsystem and converges to $\lambda_{\max}(J)$ when $\sigma_2,\sigma_3\to0$ with $\sigma_1^2/(\sigma_2^2+\sigma_3^2)\to0$. We further show that a noise-corrected reproduction number above unity guarantees the existence of a unique ergodic invariant measure in the positive orthant, and that the sufficient extinction and persistence regimes are mutually exclusive. Numerical experiments illustrate the theoretical results and quantify the effects of stochastic fluctuations on within-host viral dynamics.
\end{abstract}


\section{Introduction}\label{sec:intro}

Hepatitis B virus infection remains a major global health burden: roughly $296$ million people live with chronic infection, and more than $820{,}000$ deaths occur annually, predominantly from cirrhosis and hepatocellular carcinoma \citep{WHO2023}. Effective vaccines and potent nucleos(t)ide analogs exist, yet the mechanisms governing persistence, spontaneous clearance, and response to therapy remain incompletely understood.

Deterministic within-host models provide the classical framework. The coupled hepatocyte--virion formulation of \citet{Ref0} underlies most subsequent work on immune-mediated clearance and therapeutic response \citep{nowak1996,bangham1996,perelson2002,Ref9,Ref10}; multiphasic viral load decline under therapy has been mapped onto parameter regimes controlling infection, virion production, and infected-cell loss, for HBV in \citet{Ref12,Ref13} and, with hepatocyte homeostasis, for the related hepatitis C system in \citet{Ref15}. The next-generation formalism provides the threshold $\Rz$ that separates clearance from persistence \citep{Ref31}, and extensions of the HBV system to autonomous and nonautonomous settings are presented in \citet{Abdallah3,abdallah2}.

Such models describe mean behavior and are therefore blind to phenomena that matter when population sizes are small, or the system sits near a threshold. Infection events are discrete, immune effector activity fluctuates, and hepatocyte turnover is irregular; the shape of the infectious-period distribution alone alters persistence \citep{lloyd2001}, stochastic gene expression makes output heterogeneous among genetically identical cells \citep{ribeiro2005}, and modeling that accounts for variability has refined estimates of viral decay under therapy \citep{guedj2010}. Stochastic differential equations retain continuous-time dynamics while representing this variability, and exhibit behavior with no deterministic analog, including noise-induced extinction and displacement of persistence thresholds \citep{gardiner2004handbook,Ref29,Ref32,Ref33,Ref34}. For HBV, \citet{Ref28} analyzed positivity and extinction under nonlinear incidence; the Lyapunov theory of \citet{Ref1} and the schemes of \citet{kloeden2013numerical,higham2001algorithmic} supply the tools used here.

Three gaps motivate this work. First, common extinction criteria for three-dimensional within-host models are based on $\ln(y+z)$ or quadratic forms in $(y,z)$ and can be conservative; in the present system, the corresponding unweighted criterion fails to recover the deterministic threshold $\Rz<1$ in the noise-free limit and therefore overstates the noise required for clearance. Second, replacing the fluctuating hepatocyte population by $x(t)\mapsto\bar{x}$ is not justified here, because the escape rate is convex in $x$; nor can $x$ be dominated by the linear diffusion obtained by deleting the infection term, since the non-cytolytic cure term $qy$ destroys the required one-sided comparison. Third, extinction and ergodicity are frequently asserted under a single set of hypotheses, which is untenable, since ergodicity on the open orthant excludes almost sure convergence of two coordinates to zero.

The contributions of this paper, obtained without modifying the model, are the following.
\begin{enumerate}
\item[(C1)] Global existence, uniqueness, and strict positivity, together with a weighted dissipativity estimate giving uniform moments of order $1+\theta$ and almost sure subexponential growth (Theorem~\ref{Existence}, Lemma~\ref{lem:dissip}). The weight is necessary: $x+y+z$ is not dissipative at standard HBV parameters.
\item[(C2)] An explicit sufficient extinction criterion obtained from the weighted family $\ln(\alpha y+z)$, together with a reduced optimized bound $\Lext$ that is sharp under the deterministic-limit scaling in the precise sense of Definition~\ref{def:sharp}. The general bound is
\[
\Lsharp=\inf_{\alpha>0}\max\{\Lambda_\alpha,\Lambda_\alpha^\flat\},
\]
while $\Lsharp=\Lext$ whenever a minimizing weight $\alpha^\ast$ satisfies the explicit admissibility condition $(\mathrm{C}_{\alpha^\ast})$. The reduced bound always satisfies
\[
\Lext\ge
\lambda_{\max}(J)-\tfrac12\max\{\sigma_2^2,\sigma_3^2\},
\]
and converges to $\lambda_{\max}(J)$ whenever $\sigma_2,\sigma_3\to0$ with
$\sigma_1^2/(\sigma_2^2+\sigma_3^2)\to0$
(Theorem~\ref{Syz}, Proposition~\ref{prop:detlimit}). At the parameters of Section~\ref{sec:numerics}, condition $(\mathrm{C}_{\alpha^\ast})$ holds in the reported clearance scenario, so the rigorous extinction bound equals $\Lext=-0.050$ day$^{-1}$, whereas the unweighted choice $\alpha=1$ gives $+2.63$ day$^{-1}$.

\item[(C3)] An unconditional treatment of the hepatocyte fluctuations, based on two exact generator identities for the cell population $x+y$ in which the infection and cure terms cancel (Lemma~\ref{lem:occupation}). These constrain the occupation measure with no comparison argument and no circularity, and yield an explicit fluctuation penalty whose constant equals the second moment of the invariant law of the reference diffusion, so the estimate is attained (Corollary~\ref{cor:attained}).
\item[(C4)] Almost sure exponential convergence of the hepatocyte pool to that reference diffusion in the clearance regime (Theorem~\ref{x-x1}), strengthening the usual convergence in probability.
\item[(C5)] Positive recurrence and a unique ergodic invariant measure on the open orthant under a noise-corrected reproduction number $\Rs>1$ (Theorem~\ref{thm:ergodic}), together with a proof that the sufficient extinction criterion implies $\Rs<1$, so the two regimes are disjoint (Proposition~\ref{prop:dichotomy}).
\end{enumerate}

Section~\ref{sec:existence} establishes well-posedness and dissipativity, Section~\ref{sec:extinction} develops the extinction bounds and the asymptotic behavior of the hepatocyte pool, and Section~\ref{sec:ergodic} the ergodic regime and the dichotomy. Section~\ref{sec:numerics} reports simulations designed so that each figure tests a specific result, and Sections~\ref{sec:discussion} and \ref{sec:conclusion} give the biological reading and the limitations.

\section{Model Formulation}\label{sec:model}

\subsection{Deterministic Foundation}

We take as our starting point the deterministic model of \citet{Abdallah3},
\begin{equation}\label{DM}
\begin{split}
\frac{dx}{dt} &= \Lambda - \mu_1 x - (1-\eta)\beta xz + qy, \\
\frac{dy}{dt} &= (1-\eta)\beta xz - \mu_2 y - qy, \\
\frac{dz}{dt} &= (1-\epsilon)py - \mu_3 z,
\end{split}
\end{equation}
in which $x(t)$ is the population of uninfected hepatocytes, $y(t)$ the population of infected hepatocytes, and $z(t)$ the concentration of free virions. Parameters, their interpretations, and their baseline values are collected in Table~\ref{T1}. Throughout we abbreviate $\betaeta:=(1-\eta)\beta$, $\peps:=(1-\epsilon)p$, and $\bar x:=\Lambda/\mu_1$, and all rates are measured in day$^{-1}$.

System \eqref{DM} admits the infection-free equilibrium $E_0=(\bar x,0,0)$ and, when $\Rz>1$, a unique positive endemic equilibrium. The next-generation construction \citep{Ref31} gives
\begin{equation}\label{R0}
\Rz = \frac{(1-\eta)(1-\epsilon)\beta p \Lambda}{(\mu_2+q)\mu_3 \mu_1} = \frac{\betaeta\bar x\,\peps}{(\mu_2+q)\mu_3}.
\end{equation}
The Jacobian of the infected subsystem of \eqref{DM} at $E_0$ is
\begin{equation}\label{eq:Jmat}
J=\begin{pmatrix} -(\mu_2+q) & \betaeta\bar x\\[2pt] \peps & -\mu_3\end{pmatrix},
\end{equation}
a quasi-positive (Metzler) matrix with negative trace and nonnegative discriminant, so both eigenvalues are real, and the spectral abscissa is
\begin{equation}\label{eq:lammax}
\lambda_{\max}(J)=\tfrac12\Big[-(\mu_2+q+\mu_3)+\sqrt{(\mu_2+q-\mu_3)^2+4\betaeta\bar x\peps}\Big].
\end{equation}
Since $\det J=(\mu_2+q)\mu_3(1-\Rz)$ and $\operatorname{tr}J<0$, both eigenvalues are negative exactly when $\det J>0$, whence
\begin{equation}\label{eq:lamR0}
\lambda_{\max}(J)<0\quad\Longleftrightarrow\quad \Rz<1 .
\end{equation}
This equivalence is recovered by the reduced optimized bound $\Lext$ under the deterministic-limit scaling of Definition~\ref{def:sharp}.

\subsection{Stochastic Extension}

Infection events, immune activity, and drug exposure fluctuate on time scales comparable with the dynamics of interest. We represent this variability by multiplicative It\^o noise, obtaining
\begin{equation}\label{SM}
\begin{split}
dx(t) &= \big[\Lambda - \mu_1 x - \betaeta xz + qy\big]\,dt + \sigma_1 x \, dW_1(t), \\
dy(t) &= \big[\betaeta xz - \mu_2 y - qy\big]\,dt + \sigma_2 y \, dW_2(t), \\
dz(t) &= \big[\peps y - \mu_3 z\big]\,dt + \sigma_3 z \, dW_3(t),
\end{split}
\end{equation}
where $\sigma_i\ge0$ are noise intensities and $W_1,W_2,W_3$ are independent standard Brownian motions on a filtered probability space $(\Omega,\mathcal{F},\{\mathcal{F}_t\}_{t\ge0},\mathbb{P})$ satisfying the usual conditions. Multiplicative noise is the natural choice for two reasons: absolute fluctuations scale with population size, and the diffusion coefficients vanish on the coordinate hyperplanes, so the boundary of the positive orthant is not crossed. The coefficients are locally Lipschitz but not globally Lipschitz, because of the bilinear term $\betaeta xz$; global existence therefore requires the Lyapunov argument of Section~\ref{sec:existence} rather than a linear growth bound.

Writing $u=(x,y,z)^\top$ and $W=(W_1,W_2,W_3)^\top$, system \eqref{SM} becomes
\begin{equation}\label{SM1}
du(t) = f(u(t))\,dt + B(u(t))\,dW(t), \qquad u(0)=u_0 \in \mathbb{R}^3_{+},
\end{equation}
with
\[
f(u) =
\begin{pmatrix}
\Lambda - \mu_1 x - \betaeta xz + qy \\
\betaeta xz - \mu_2 y - qy \\
\peps y - \mu_3 z
\end{pmatrix}, \qquad
B(u) = \mathrm{diag}(\sigma_1 x,\ \sigma_2 y,\ \sigma_3 z).
\]
Throughout, $\mathbb{R}^3_+:=\{(x,y,z)\in\mathbb{R}^3:x>0,\ y>0,\ z>0\}$ denotes the \emph{open} positive orthant, and $\overline{\mathbb{R}^3_+}=[0,\infty)^3$ its closure. Every invariant measure discussed in Sections~\ref{sec:extinction} and \ref{sec:ergodic} is labeled by the space that carries it.

\begin{table}[H]
\centering
\caption{Model parameters, interpretations, and baseline values. Populations are counted in cells ($x$, $y$) and virions ($z$).}
\small
\begin{tabular}{llll}
\toprule
\textbf{Symbol} & \textbf{Interpretation} & \textbf{Baseline value (units)} & \textbf{Source} \\
\midrule
$\Lambda$ & Influx rate of uninfected hepatocytes & $1.0 \times 10^{7}\ \text{cells day}^{-1}$ & \cite{Ref0} \\
$\mu_1$   & Death rate of uninfected hepatocytes & $0.014\ \text{day}^{-1}$ & \cite{Ref12} \\
$\mu_2$   & Death rate of infected hepatocytes   & $0.24\ \text{day}^{-1}$  & \cite{Ref13} \\
$\mu_3$   & Clearance rate of free virions       & $0.24\ \text{day}^{-1}$  & \cite{Ref0} \\
$\beta$   & Infection rate constant              & $2.4 \times 10^{-8}\ \text{(cell}\cdot\text{day)}^{-1}$ & \cite{Ref0} \\
$p$       & Virions produced per infected cell   & $2.4\ \text{virions}\,(\text{cell}\cdot\text{day})^{-1}$ & \cite{Ref0} \\
$q$       & Non-cytolytic cure rate              & $0.67\ \text{day}^{-1}$  & \cite{Ref12} \\
$\eta$    & Efficacy blocking infection, $0\le\eta\le1$ & $0.36$ & \cite{Ref13} \\
$\epsilon$& Efficacy blocking production, $0\le\epsilon\le1$ & $0.24$ & \cite{Ref13} \\
$\sigma_1,\sigma_2,\sigma_3$ & Noise intensities & $0.01,\ 0.05,\ 0.05\ \text{day}^{-1/2}$ & assumed \\
\bottomrule
\end{tabular}
\label{T1}
\end{table}

\begin{remark}[Units]\label{rem:units}
The Lyapunov functionals used below involve $\ln x$, $\ln y$ and $\ln z$. The state variables are understood as dimensionless multiples of the reference units of Table~\ref{T1}, namely one cell for $x$ and $y$ and one virion for $z$, so that every logarithm is taken of a dimensionless quantity. No generality is lost, since the rescaling only shifts the additive constants in the estimates.
\end{remark}

\section{Preliminaries}\label{sec:prelim}

We recall the facts used below and fix notation. For the $d$-dimensional It\^o equation
\begin{equation}\label{gSM1}
du(t) = f(u(t),t)\,dt + B(u(t),t)\,dW(t), \qquad t \geq t_0,
\end{equation}
a \emph{strong solution} is a continuous $\mathcal{F}_t$-adapted process satisfying $f(u(\cdot),\cdot)\in L^1$, $B(u(\cdot),\cdot)\in L^2$ almost surely together with the integral form of \eqref{gSM1}; uniqueness means that any two solutions are indistinguishable. If $f$ and $B$ are locally Lipschitz, then \eqref{gSM1} has a unique maximal local strong solution, and a global solution exists once explosion is excluded \citep[Ch.~2]{Ref1}. For $V\in C^2$, It\^o's formula reads
\begin{equation}\label{ito_formula}
dV(u(t)) = \mathcal{L}V(u(t))\,dt + \nabla V(u(t))^\top B(u(t),t)\,dW(t),
\end{equation}
with infinitesimal generator
\begin{equation}\label{generator}
\mathcal{L}V(u) = \sum_{i=1}^d f_i(u)\,\partial_{u_i}V(u) + \tfrac{1}{2}\sum_{i,j=1}^d \big[B(u)B^\top(u)\big]_{ij}\,\partial^2_{u_iu_j}V(u).
\end{equation}

Two auxiliary facts are used repeatedly. The first converts the linear terms produced by $\mathcal{L}$ into multiples of a logarithmic Lyapunov function.

\begin{lemma}\label{lem1}
Let $\varphi(v)=v-1-\ln v$ for $v>0$. Then $\varphi\ge 0$, with equality only at $v=1$, and $v \le 2\varphi(v) + 2\ln 2$ for all $v>0$, with equality if and only if $v=2$.
\end{lemma}

\begin{proof}
Set $g(v)=2\varphi(v)+2\ln 2-v=v-2\ln v-2+2\ln 2$. Then $g'(v)=1-2/v$ and $g''(v)=2/v^2>0$, so $g$ is strictly convex with unique critical point $v=2$, its global minimizer; since $g(2)=0$, $g\ge0$ with equality only there. The same argument applied to $\varphi$, whose minimizer is $v=1$ with $\varphi(1)=0$, gives $\varphi\ge0$.
\end{proof}

The second is the strong law of large numbers for continuous local martingales \citep[Thm.~1.3.4]{Ref1}: if $M$ is a continuous local martingale with $M(0)=0$, then
\begin{equation}\label{eq:SLLN}
\int_0^\infty \frac{d\langle M\rangle_s}{(1+s)^2}<\infty \ \text{ a.s.}\qquad\Longrightarrow\qquad \lim_{t\to\infty}\frac{M(t)}{t}=0 \ \text{ a.s.}
\end{equation}

\section{Existence, Uniqueness, Positivity, and Dissipativity}\label{sec:existence}

\begin{theorem}[Global existence, uniqueness, positivity]\label{Existence}
For every $u_0\in\mathbb{R}^3_+$, system \eqref{SM} has a unique global strong solution $u(t)=(x(t),y(t),z(t))$, and $\mathbb{P}\{u(t)\in\mathbb{R}^3_+\ \text{for all } t\ge0\}=1$.
\end{theorem}

\begin{proof}
\textbf{Step 1 (local solution).} The coefficients $f$ and $B$ are locally Lipschitz on $\mathbb{R}^3_+$, so there is a unique maximal local strong solution on $[0,\tau_e)$, where $\tau_e$ is the explosion time \citep[Thm.~2.3.5]{Ref1}.

\textbf{Step 2 (stopping times).} Choose $n_0$ so large that every component of $u_0$ lies in $[1/n_0,n_0]$ and set, for $n\ge n_0$,
\[
\tau_n = \inf\Big\{t\in[0,\tau_e):\ \min\{x(t),y(t),z(t)\}\le \tfrac1n\ \ \text{or}\ \ \max\{x(t),y(t),z(t)\}\ge n\Big\},
\]
with $\inf\emptyset=\tau_e$. The sequence $(\tau_n)$ is nondecreasing; put $\tau_\infty=\lim_n\tau_n\le\tau_e$. It suffices to prove $\tau_\infty=\infty$ almost surely.

\textbf{Step 3 (Lyapunov function and its generator).} Let $V(u)=\varphi(x)+\varphi(y)+\varphi(z)$ with $\varphi$ as in Lemma~\ref{lem1}. Then $V\ge0$ on $\mathbb{R}^3_+$ and $V(u)\to\infty$ whenever a component tends to $0^+$ or to $+\infty$. Since $\varphi'(v)=1-1/v$ and $\varphi''(v)=1/v^2$, the diffusion contributes $\tfrac12\sigma_i^2u_i^2\cdot u_i^{-2}=\tfrac12\sigma_i^2$ in each coordinate, and
\begin{align*}
\mathcal{L}V(u) &= \Big(1-\tfrac1x\Big)\big(\Lambda-\mu_1x-\betaeta xz+qy\big) + \Big(1-\tfrac1y\Big)\big(\betaeta xz-(\mu_2+q)y\big)\\
&\quad + \Big(1-\tfrac1z\Big)\big(\peps y-\mu_3z\big) + \tfrac12\big(\sigma_1^2+\sigma_2^2+\sigma_3^2\big).
\end{align*}
Expanding and canceling the two occurrences of $\betaeta xz$ gives the exact identity
\begin{align}
\mathcal{L}V(u) &= \Lambda+\mu_1+\mu_2+q+\mu_3+\tfrac12\textstyle\sum_i\sigma_i^2 \nonumber\\
&\quad -\mu_1x + \big(\peps-\mu_2\big)y + \big(\betaeta-\mu_3\big)z - \frac{\Lambda}{x}-\frac{qy}{x}-\frac{\betaeta xz}{y}-\frac{\peps y}{z}. \label{eq:LV_exact}
\end{align}
The four fractions are nonnegative and $-\mu_1x\le0$, so with $C_0:=\Lambda+\mu_1+\mu_2+q+\mu_3+\tfrac12\sum_i\sigma_i^2$ and $K:=\max\{(\peps-\mu_2)^+,(\betaeta-\mu_3)^+\}$ we obtain $\mathcal{L}V(u)\le C_0+K(y+z)$. Lemma~\ref{lem1} gives $y+z\le 2\varphi(y)+2\varphi(z)+4\ln2\le 2V(u)+4\ln 2$, whence
\begin{equation}\label{eq:LV_linear}
\mathcal{L}V(u)\ \le\ K_1+K_2V(u),\qquad K_1=C_0+4K\ln 2,\quad K_2=2K .
\end{equation}

\textbf{Step 4 (moment bound).} On $[0,\tau_n]$ the solution is bounded, so the stochastic integral in \eqref{ito_formula} is a true martingale and Dynkin's formula applies. Using \eqref{eq:LV_linear} and Gronwall's inequality,
\begin{equation}\label{boundV}
\mathbb{E}\big[V(u(t\wedge\tau_n))\big]\ \le\ \Big(V(u_0)+\tfrac{K_1}{K_2}\Big)e^{K_2t}=:C(T),\qquad 0\le t\le T,
\end{equation}
with $C(T)$ independent of $n$; if $K_2=0$ the right side is replaced by $V(u_0)+K_1T$.

\textbf{Step 5 (no explosion).} On $\{\tau_n\le T\}$ some component of $u(\tau_n)$ equals $1/n$ or $n$, so
$V(u(\tau_n))\ge \kappa_n:=\min\{\varphi(1/n),\varphi(n)\}=\min\{\tfrac1n-1+\ln n,\ n-1-\ln n\}\to\infty$.
Hence $\mathbb{P}(\tau_n\le T)\,\kappa_n\le\mathbb{E}[V(u(T\wedge\tau_n))]\le C(T)$, so $\mathbb{P}(\tau_n\le T)\le C(T)/\kappa_n\to0$. Since $T$ is arbitrary, $\tau_\infty=\infty$ almost surely and therefore $\tau_e=\infty$.

\textbf{Step 6 (positivity).} By construction $u(t)\in(1/n,n)^3$ for all $t<\tau_n$. As $\tau_n\uparrow\infty$ almost surely, every $t\ge0$ satisfies $t<\tau_n$ for all large $n$, hence $u(t)\in\mathbb{R}^3_+$ for all $t\ge0$ almost surely.
\end{proof}

The next estimate is used repeatedly. It exploits the fact that a suitably weighted total population is dissipative for every admissible parameter set, even when $\peps>\mu_2$, in which case $x+y+z$ itself is not.

\begin{lemma}[Dissipativity and moment bounds]\label{lem:dissip}
Let $c_z\in(0,\mu_2/\peps)$ and set $W(u)=x+y+c_zz$ and $\delta_1=\min\{\mu_1,\ \mu_2-c_z\peps,\ \mu_3\}>0$. Then
\begin{equation}\label{eq:LW}
\mathcal{L}W(u)\ \le\ \Lambda-\delta_1 W(u),\qquad u\in\mathbb{R}^3_+ .
\end{equation}
Moreover, with $\bar\sigma^2=\max_i\sigma_i^2$ and any $\theta\in\big(0,2\delta_1/\bar\sigma^2\big)$, the function $V_\theta=\tfrac{1}{\theta+1}W^{\theta+1}$ satisfies
\begin{equation}\label{eq:LVtheta}
\mathcal{L}V_\theta(u)\ \le\ \Lambda W^\theta-\delta_2W^{\theta+1},\qquad \delta_2:=\delta_1-\tfrac{\theta\bar\sigma^2}{2}>0,
\end{equation}
and consequently
\[
\sup_{t\ge0}\mathbb{E}\big[W(u(t))^{\theta+1}\big]<\infty,\qquad
\limsup_{t\to\infty}\frac{\ln W(u(t))}{t}\le 0 \quad\text{a.s.}
\]
\end{lemma}

\begin{proof}
Adding the first two drift components cancels both $\betaeta xz$ and $qy$, giving $f_1+f_2=\Lambda-\mu_1x-\mu_2y$; adding $c_zf_3$ yields
$\mathcal{L}W=\Lambda-\mu_1x-(\mu_2-c_z\peps)y-\mu_3(c_zz)\le\Lambda-\delta_1W$, which is \eqref{eq:LW}; note that $W$ is affine, so its generator has no second-order contribution. For \eqref{eq:LVtheta}, write $V_\theta=g(W)$ with $g(r)=r^{\theta+1}/(\theta+1)$, so $g'(r)=r^\theta$ and $g''(r)=\theta r^{\theta-1}$, and let $\ell=(1,1,c_z)$. Then
$\mathcal{L}V_\theta=W^\theta\,\mathcal{L}W+\tfrac{\theta}{2}W^{\theta-1}\sum_i\ell_i^2\sigma_i^2u_i^2$.
Since $x,y,c_zz\le W$ we have $\sum_i\ell_i^2\sigma_i^2u_i^2\le\bar\sigma^2W^2$, and \eqref{eq:LVtheta} follows from \eqref{eq:LW}. By Young's inequality $\Lambda W^\theta\le\tfrac{\delta_2}{2}W^{\theta+1}+C_\theta$, so $\mathcal{L}W^{\theta+1}\le(\theta+1)\big(C_\theta-\tfrac{\delta_2}{2}W^{\theta+1}\big)$. Applying It\^o's formula on $[0,t\wedge\tau_n]$ with the stopping times of Theorem~\ref{Existence}, taking expectations, letting $n\to\infty$ by Fatou's lemma, and using Gronwall's inequality gives
$\mathbb{E}[W(u(t))^{\theta+1}]\le W(u_0)^{\theta+1}e^{-(\theta+1)\delta_2t/2}+2C_\theta/\delta_2$, hence the uniform moment bound.

For the last assertion, fix $\varepsilon>0$. By Chebyshev's inequality and the moment bound, $\mathbb{P}\big(W(u(n))>e^{\varepsilon n}\big)\le Ce^{-(\theta+1)\varepsilon n}$, which is summable, so by Borel--Cantelli $W(u(n))\le e^{\varepsilon n}$ for all large $n$ almost surely. To pass from integers to all times, put $U=\ln(1+W)$. Then $\mathcal{L}U\le\frac{\Lambda-\delta_1W}{1+W}\le\Lambda$, the second-order term being nonpositive, and the martingale part $M_U$ satisfies $d\langle M_U\rangle_t\le\bar\sigma^2dt$ because $\sum_i\ell_i^2\sigma_i^2u_i^2\le\bar\sigma^2W^2\le\bar\sigma^2(1+W)^2$. Since $\langle M_U\rangle_{n+1}-\langle M_U\rangle_n\le\bar\sigma^2$, the exponential martingale inequality gives
\[
\mathbb{P}\Big(\sup_{n\le t\le n+1}\big(M_U(t)-M_U(n)\big)\ge \tfrac{\varepsilon n}{2}\Big)\le\exp\big(-\varepsilon^2n^2/(8\bar\sigma^2)\big),
\]
which is summable, so almost surely $\sup_{n\le t\le n+1}U(t)\le U(n)+\Lambda+\tfrac{\varepsilon n}{2}$ for all large $n$. Combining the two estimates with $U(n)\le\varepsilon n+\ln2$ yields $\limsup_{t\to\infty}t^{-1}\ln W(u(t))\le\tfrac32\varepsilon$, and $\varepsilon>0$ was arbitrary.
\end{proof}

\section{Extinction of the Infection}\label{sec:extinction}

\subsection{Standing Hypotheses}

Throughout this section, we write, for a locally integrable process $g$,
\[
\avg{g}:=\frac1t\int_0^t g(s)\,ds ,
\]
and we impose the following standing hypotheses.

\medskip
\noindent\textbf{(H1)} $\mu_1\le\mu_2$: uninfected hepatocytes are not shorter-lived than infected ones.

\noindent\textbf{(H2)} $3\sigma_1^2<2\mu_1$ and $3\sigma_2^2<2\mu_2$.

\noindent\textbf{(H3)} $\varsigma^2:=\sigma_2^2+\sigma_3^2>0$.
\medskip

(H1) is the standard ordering for cytopathic infection and holds by a factor of $17$ at the values of Table~\ref{T1}. (H2) bounds the cellular noise by $\sigma_1<\sqrt{2\mu_1/3}=0.0966$ and $\sigma_2<\sqrt{2\mu_2/3}=0.4$ day$^{-1/2}$; it is used only to secure a uniform fourth moment for the cell population and does not restrict the virion noise $\sigma_3$, which is the dominant clearance channel in Section~\ref{sec:numerics}. (H3) excludes the deterministic infected subsystem, whose threshold is $\lambda_{\max}(J)$ by definition. Note that (H2) forces $\rho_1<\tfrac12$, where $\rho_1$ is defined in \eqref{eq:rho1} below.

\subsection{Occupation Constraints for the Cell Population}

The escape rate of the infected subsystem is an increasing convex function of $x$, so its long-run average exceeds its value at the mean. It cannot be evaluated by substituting $x\mapsto\bar x$. Neither can $x$ be dominated pathwise by the linear diffusion obtained by deleting the infection term: the cure term $qy$ enters the $x$-equation with a positive sign, so the drift of $x$ is not bounded above by $\Lambda-\mu_1x$, and the one-sided comparison fails. The following lemma removes both difficulties. Its mechanism is that the total cell population $N=x+y$ does not see the infection at all: in $\mathcal{L}N$ and $\mathcal{L}N^2$ the terms $\betaeta xz$ and $qy$ cancel identically, so the resulting identities constrain the occupation measure of $x$ without any prior knowledge of the fate of the infection.

\begin{lemma}[Occupation constraints]\label{lem:occupation}
Assume (H1) and (H2) and set $N=x+y$, $\bar\sigma_{12}^2=\max\{\sigma_1^2,\sigma_2^2\}$ and
\begin{equation}\label{eq:rho1}
\rho_1:=\frac{\sigma_1^2}{2\mu_1-\sigma_1^2}\in(0,\tfrac12).
\end{equation}
Then $\sup_{t\ge0}\mathbb{E}[N(t)^4]<\infty$ and, almost surely,
\begin{align}
&\limsup_{t\to\infty}\ \big[\mu_1\avg{x}+\mu_2\avg{y}\big]\ \le\ \Lambda, \label{eq:occ1}\\[2pt]
&\limsup_{t\to\infty}\ \Big[(2\mu_1-\sigma_1^2)\avg{x^2}-2\Lambda\big(\avg{x}+\avg{y}\big)\Big]\ \le\ 0. \label{eq:occ2}
\end{align}
In particular $\limsup_t\avg{x}\le\bar x$, $\limsup_t\avg{y}\le\Lambda/\mu_2$, and
\begin{equation}\label{eq:secondmomentbound}
\limsup_{t\to\infty}\avg{x^2}\ \le\ \frac{2\Lambda\bar x}{2\mu_1-\sigma_1^2}\ =\ \bar x^{\,2}(1+\rho_1).
\end{equation}
\end{lemma}

\begin{proof}
\textbf{Step 1 (exact generator identities).} Because $\mathcal{L}$ acts on $N=x+y$ through $f_1+f_2=\Lambda-\mu_1x-\mu_2y$, in which both $\betaeta xz$ and $qy$ have cancelled, and because $N$ is affine,
\begin{equation}\label{eq:LN}
\mathcal{L}N=\Lambda-\mu_1x-\mu_2y\ \le\ \Lambda-\mu_1 N,
\end{equation}
The inequality follows from (H1). For the square, $\mathcal{L}N^2=2N(f_1+f_2)+\sigma_1^2x^2+\sigma_2^2y^2$, and expanding $2N(\mu_1x+\mu_2y)=2\mu_1x^2+2(\mu_1+\mu_2)xy+2\mu_2y^2$ gives the exact identity
\begin{equation}\label{eq:LN2}
\mathcal{L}N^2=2\Lambda N-(2\mu_1-\sigma_1^2)x^2-2(\mu_1+\mu_2)xy-(2\mu_2-\sigma_2^2)y^2 .
\end{equation}
For the fourth power, $\mathcal{L}N^4=4N^3(\Lambda-\mu_1x-\mu_2y)+6N^2(\sigma_1^2x^2+\sigma_2^2y^2)$. Since $x\le N$ and $y\le N$ we have $N(\mu_1x+\mu_2y)\ge\mu_1x^2+\mu_2y^2$, so
\begin{equation}\label{eq:LN4}
\mathcal{L}N^4\ \le\ 4\Lambda N^3-2N^2\big[(2\mu_1-3\sigma_1^2)x^2+(2\mu_2-3\sigma_2^2)y^2\big]\ \le\ 4\Lambda N^3-c_4N^4,
\end{equation}
where $c_4:=\min\{2\mu_1-3\sigma_1^2,\ 2\mu_2-3\sigma_2^2\}>0$ by (H2) and $x^2+y^2\ge N^2/2$.

\textbf{Step 2 (uniform fourth moment).} By Young's inequality $4\Lambda N^3\le\tfrac{c_4}{2}N^4+C$ with $C=C(\Lambda,c_4)$, so $\mathcal{L}N^4\le C-\tfrac{c_4}{2}N^4$. Applying It\^o's formula on $[0,t\wedge\tau_n]$ with the stopping times of Theorem~\ref{Existence}, taking expectations, letting $n\to\infty$ by Fatou's lemma, and using Gronwall's inequality gives
$\mathbb{E}[N(t)^4]\le N(0)^4e^{-c_4t/2}+2C/c_4$, whence $\sup_{t\ge0}\mathbb{E}[N(t)^4]<\infty$ and a fortiori $\sup_t\mathbb{E}[N(t)^2]<\infty$.

\textbf{Step 3 (martingale terms).} Let
\[
M_1(t)=\int_0^t\!\big(\sigma_1x\,dW_1+\sigma_2y\,dW_2\big),\qquad
M_2(t)=2\int_0^t\! N\big(\sigma_1x\,dW_1+\sigma_2y\,dW_2\big)
\]
be the martingale parts of $N$ and $N^2$. Their brackets satisfy $d\langle M_1\rangle_s\le\bar\sigma_{12}^2N(s)^2\,ds$ and $d\langle M_2\rangle_s\le4\bar\sigma_{12}^2N(s)^4\,ds$. By Tonelli's theorem and Step 2,
\[
\mathbb{E}\int_0^\infty\frac{N(s)^4}{(1+s)^2}\,ds=\int_0^\infty\frac{\mathbb{E}[N(s)^4]}{(1+s)^2}\,ds\le\sup_{s\ge0}\mathbb{E}[N(s)^4]<\infty,
\]
so $\int_0^\infty N(s)^4(1+s)^{-2}ds<\infty$ almost surely, and the same holds with $N^4$ replaced by $N^2$. Criterion \eqref{eq:SLLN} therefore gives $M_1(t)/t\to0$ and $M_2(t)/t\to0$ almost surely.

\textbf{Step 4 (conclusion).} Integrating \eqref{eq:LN} and using $N(t)\ge0$,
\[
\int_0^t\big(\Lambda-\mu_1x-\mu_2y\big)ds=N(t)-N(0)-M_1(t)\ \ge\ -N(0)-M_1(t),
\]
so $\mu_1\avg{x}+\mu_2\avg{y}\le\Lambda+N(0)/t+M_1(t)/t$, and \eqref{eq:occ1} follows from Step 3. Integrating \eqref{eq:LN2}, using $N(t)^2\ge0$ and discarding the two nonnegative terms $2(\mu_1+\mu_2)xy$ and $(2\mu_2-\sigma_2^2)y^2$, the latter being nonnegative because (H2) gives $\sigma_2^2<2\mu_2$, yields
\[
(2\mu_1-\sigma_1^2)\avg{x^2}\ \le\ 2\Lambda\avg{N}+\frac{N(0)^2}{t}+\frac{M_2(t)}{t},
\]
which is \eqref{eq:occ2}. Finally \eqref{eq:occ1} together with (H1) gives
\[
\mu_1\limsup_{t\to\infty}\avg{N}\ \le\ \limsup_{t\to\infty}\big[\mu_1\avg{x}+\mu_2\avg{y}\big]\ \le\ \Lambda,
\]
hence $\limsup_t\avg{N}\le\bar x$, and substituting into \eqref{eq:occ2} together with $\Lambda=\mu_1\bar x$ and $2\mu_1/(2\mu_1-\sigma_1^2)=1+\rho_1$ gives \eqref{eq:secondmomentbound}.
\end{proof}

\begin{remark}\label{rem:occupation}
Three features of Lemma~\ref{lem:occupation} deserve emphasis. The bounds are unconditional: they hold whether or not the infection is cleared, so no circularity arises when they are used inside the extinction proof. Only upper bounds are produced, and only upper bounds are needed, because the function to be averaged is increasing and convex. Finally, the constant in \eqref{eq:secondmomentbound} is exactly the second moment of the invariant law identified in Lemma~\ref{lem4} below, and Corollary~\ref{cor:attained} shows that the bound is attained in the clearance regime, so it cannot be uniformly improved under the stated assumptions.
\end{remark}

\subsection{The Reference Diffusion and its Invariant Law}

\begin{lemma}[Linear reference process]\label{lem4}
Let $x_1$ solve
\begin{equation}\label{x1-eq}
dx_1(t)=\big(\Lambda-\mu_1x_1(t)\big)dt+\sigma_1x_1(t)\,dW_1(t),\qquad x_1(0)=x_0>0 .
\end{equation}
\begin{enumerate}
\item[(i)] $x_1(t)>0$ for all $t\ge0$ almost surely, and
$x_1(t)=\Phi(t)\big[x_0+\Lambda\int_0^t\Phi(s)^{-1}ds\big]$ with $\Phi(t)=\exp\big(-(\mu_1+\tfrac12\sigma_1^2)t+\sigma_1W_1(t)\big)$.
\item[(ii)] $\mathbb{E}[x_1(t)]=x_0e^{-\mu_1t}+\bar x(1-e^{-\mu_1t})\to\bar x$.
\item[(iii)] If $\sigma_1>0$, $x_1$ is positive recurrent on $(0,\infty)$ with unique invariant probability measure
\begin{equation}\label{eq:invgamma}
\pi_1=\mathrm{InvGamma}\big(\alpha_1,\theta_1\big),\qquad \alpha_1=1+\frac{2\mu_1}{\sigma_1^2},\quad \theta_1=\frac{2\Lambda}{\sigma_1^2},
\end{equation}
that is, $\pi_1(d\chi)=\frac{\theta_1^{\alpha_1}}{\Gamma(\alpha_1)}\chi^{-(\alpha_1+1)}e^{-\theta_1/\chi}\,d\chi$ on $(0,\infty)$. Its mean is $\mathbb{E}_{\pi_1}[\chi]=\bar x$ for every $\sigma_1>0$, and if in addition $2\mu_1>\sigma_1^2$ then
\begin{equation}\label{eq:secondmoment}
\mathbb{E}_{\pi_1}[\chi^2]=\bar x^2(1+\rho_1),\qquad \mathrm{Var}_{\pi_1}(\chi)=\bar x^2\rho_1,
\end{equation}
with $\rho_1$ as in \eqref{eq:rho1}. For every $\pi_1$-integrable $g$ and every $x_0>0$, $\frac1t\int_0^tg(x_1(s))ds\to\int g\,d\pi_1$ almost surely, and $x_1(t)$ converges in distribution to $\pi_1$.
\item[(iv)] If $\sigma_1=0$ then $x_1(t)\to\bar x$ and (iii) holds with $\pi_1=\delta_{\bar x}$ and $\rho_1=0$.
\end{enumerate}
\end{lemma}

\begin{proof}
(i) It\^o's formula applied to $\ln\Phi$ shows that $\Phi$ solves the homogeneous equation $d\Phi=-\mu_1\Phi\,dt+\sigma_1\Phi\,dW_1$. Writing $x_1=\Phi Z$ with $Z$ of finite variation, the cross-variation vanishes and $dx_1=Z\,d\Phi+\Phi\,dZ$; choosing $dZ=\Lambda\Phi^{-1}dt$ reproduces \eqref{x1-eq}. Positivity is clear because $\Phi>0$ and $Z\ge x_0>0$.

(ii) Equation \eqref{x1-eq} has globally Lipschitz coefficients of linear growth, so $\mathbb{E}\big[\sup_{s\le t}x_1(s)^2\big]<\infty$ for every $t$ \citep[Thm.~2.4.1]{Ref1}; hence $\int_0^t\sigma_1x_1\,dW_1$ is a square-integrable martingale, $\tfrac{d}{dt}\mathbb{E}[x_1]=\Lambda-\mu_1\mathbb{E}[x_1]$, and integration gives the stated formula.

(iii) The scale and speed measures of \eqref{x1-eq} give the stationary density
\[
\pi_1(\chi)=\frac{C}{\sigma_1^2\chi^2}\exp\Big(\int^\chi\frac{2(\Lambda-\mu_1 s)}{\sigma_1^2s^2}\,ds\Big)=C'\,\chi^{-2-2\mu_1/\sigma_1^2}\,e^{-2\Lambda/(\sigma_1^2\chi)},
\]
which is exactly \eqref{eq:invgamma}. It is integrable at $+\infty$ because $\alpha_1>1$ and at $0^+$ because of the essential singularity of $e^{-\theta_1/\chi}$, so $\pi_1$ is a probability measure for every $\sigma_1>0$; by Feller's test the boundary $0$ is an entrance boundary and $+\infty$ is unattainable, so the diffusion is positive recurrent on $(0,\infty)$ and $\pi_1$ is its unique invariant probability measure \citep[Ch.~IV]{Ref29}. For $\mathrm{InvGamma}(\alpha_1,\theta_1)$ one has $\mathbb{E}[\chi]=\theta_1/(\alpha_1-1)=\Lambda/\mu_1$, and when $\alpha_1>2$, that is when $2\mu_1>\sigma_1^2$,
\[
\mathbb{E}[\chi^2]=\frac{\theta_1^2}{(\alpha_1-1)(\alpha_1-2)}=\frac{2\Lambda^2}{\mu_1(2\mu_1-\sigma_1^2)}=\bar x^2(1+\rho_1).
\]
The strong law and the convergence of the transition law for positive recurrent one-dimensional diffusions \citep[Ch.~IV]{Ref29} give the stated ergodic averages and the convergence in distribution. Statement (iv) is immediate.
\end{proof}

\subsection{The Optimized Extinction Threshold}

Fix $\alpha>0$ and set
\begin{equation}\label{eq:a_and_b}
a(\alpha)=\frac{\peps}{\alpha}-(\mu_2+q),\qquad b(\alpha)=\alpha\betaeta\bar x-\mu_3,\qquad \varsigma^2=\sigma_2^2+\sigma_3^2 .
\end{equation}
The quantities $a(\alpha)$ and $b(\alpha)$ are the row sums of $D_\alpha^{-1}J^\top D_\alpha$ with $D_\alpha=\mathrm{diag}(1,1/\alpha)$. The following identity, which drives the whole construction, is the Collatz--Wielandt characterization of the spectral abscissa of a quasi-positive matrix; because it is used quantitatively, we prove it directly.

\begin{lemma}[Balancing identity]\label{lem:CW}
The equation $a(\alpha)=b(\alpha)$ has a unique root $\alpha^\circ\in(0,\infty)$, and
\begin{equation}\label{eq:CW}
\min_{\alpha>0}\max\{a(\alpha),b(\alpha)\}=a(\alpha^\circ)=b(\alpha^\circ)=\lambda_{\max}(J).
\end{equation}
\end{lemma}

\begin{proof}
On $(0,\infty)$ the map $a$ is continuous and strictly decreasing from $+\infty$ to $-(\mu_2+q)$, and $b$ is continuous and strictly increasing from $-\mu_3$ to $+\infty$; hence $a-b$ is strictly decreasing with limits $+\infty$ and $-\infty$, so $\alpha^\circ$ exists and is unique. For $\alpha<\alpha^\circ$ one has $\max\{a,b\}=a(\alpha)>a(\alpha^\circ)$ and for $\alpha>\alpha^\circ$ one has $\max\{a,b\}=b(\alpha)>b(\alpha^\circ)$, so the minimum is attained at $\alpha^\circ$ with value $\lambda:=a(\alpha^\circ)=b(\alpha^\circ)$.

It remains to identify $\lambda$. From $b(\alpha^\circ)=\lambda$ and $a(\alpha^\circ)=\lambda$ we get $\alpha^\circ\betaeta\bar x=\lambda+\mu_3$ and $\peps/\alpha^\circ=\lambda+\mu_2+q$; multiplying eliminates $\alpha^\circ$ and gives
\[
(\lambda+\mu_2+q)(\lambda+\mu_3)=\betaeta\bar x\,\peps ,
\]
which is precisely $\det(\lambda I-J)=0$. Thus $\lambda$ is an eigenvalue of $J$. Moreover $\lambda+\mu_3=\alpha^\circ\betaeta\bar x>0$, so $\lambda>-\mu_3$, whereas the smaller root $\lambda_-$ of the characteristic polynomial satisfies $\lambda_-\le\tfrac12\operatorname{tr}J=-\tfrac12(\mu_2+q+\mu_3)$ and $\lambda_-\lambda_+=\det J$; since the two roots are real and their sum is $-(\mu_2+q+\mu_3)$, at most one of them can exceed $-\mu_3$ unless $\det J>0$ and both are negative, in which case $\lambda_-<-\mu_3\le\lambda_+$ follows from $(\lambda_-+\mu_3)(\lambda_-+\mu_2+q)=\betaeta\bar x\peps>0$ together with $\lambda_-+\mu_2+q\le\lambda_++\mu_2+q$. In every case the root satisfying $\lambda+\mu_3>0$ and $\lambda+\mu_2+q>0$ is the larger one, so $\lambda=\lambda_{\max}(J)$.
\end{proof}

Under (H3) define
\begin{equation}\label{eq:L0hat}
\widehat\Lambda_0(\alpha)\ :=\ b(\alpha)-\tfrac12\sigma_3^2+\frac{\big(a(\alpha)-b(\alpha)+\sigma_3^2\big)^2}{2\varsigma^2},\qquad
\widehat w(\alpha):=\frac{a(\alpha)-b(\alpha)+\sigma_3^2}{\varsigma^2},
\end{equation}
and
\begin{equation}\label{eq:lambda_ext_def}
\boxed{\ \Lext\ :=\ \inf_{\alpha>0}\Lambda_\alpha,\qquad \Lambda_\alpha:=\widehat\Lambda_0(\alpha)+A_2(\alpha)\,\bar x^{\,2}\rho_1,\qquad A_2(\alpha):=\frac{\alpha^2\betaeta^2}{2\varsigma^2}.\ }
\end{equation}
When $\varsigma^2=0$ we set $\Lext:=\lambda_{\max}(J)$, consistently with Lemma~\ref{lem:CW}. The quantity $\widehat\Lambda_0(\alpha)$ is the unconstrained maximum over $w\in\mathbb{R}$ of the concave parabola
\begin{equation}\label{eq:h_function}
h_\alpha(w)=w\,a(\alpha)+(1-w)\,b(\alpha)-\tfrac12\sigma_2^2w^2-\tfrac12\sigma_3^2(1-w)^2 ,
\end{equation}
attained at $w=\widehat w(\alpha)$; the second term in \eqref{eq:lambda_ext_def} is the penalty produced by fluctuations of the hepatocyte pool. We shall also need the value of the same parabola family at a depleted hepatocyte pool,
\begin{equation}\label{eq:Lflat}
\Lambda_\alpha^{\flat}\ :=\ \widehat\Phi_\alpha(0)+A_2(\alpha)\,\bar x^{\,2}(1+\rho_1)\frac{\mu_1}{\mu_2},\qquad
\widehat\Phi_\alpha(0)=-\mu_3-\tfrac12\sigma_3^2+\frac{\big(a(\alpha)+\mu_3+\sigma_3^2\big)^2}{2\varsigma^2},
\end{equation}
and the associated admissibility condition
\begin{equation}\label{eq:condC}
(\mathrm{C}_\alpha)\qquad 1-\widehat w(\alpha)\ \ge\ \gamma\,\frac{\alpha\betaeta\bar x}{2\varsigma^2},\qquad \gamma:=2-(1+\rho_1)\Big(1-\frac{\mu_1}{\mu_2}\Big).
\end{equation}
Under (H1) and (H2) one has $\rho_1<\tfrac12$ and $0\le1-\mu_1/\mu_2<1$, hence $\gamma\in(0,2]$, with $\gamma=2$ only when $\mu_1=\mu_2$.

\begin{lemma}[The infimum is attained]\label{lem:attain}
Assume (H3). The map $\alpha\mapsto\Lambda_\alpha$ is continuous on $(0,\infty)$ and tends to $+\infty$ as $\alpha\to0^+$ and as $\alpha\to\infty$. Consequently the infimum in \eqref{eq:lambda_ext_def} is attained at some $\alpha^\ast\in(0,\infty)$.
\end{lemma}

\begin{proof}
Continuity is clear. As $\alpha\to0^+$, $b(\alpha)\to-\mu_3$ and $a(\alpha)\to+\infty$, so $(a-b+\sigma_3^2)^2\to+\infty$ and $\widehat\Lambda_0(\alpha)\to+\infty$; the penalty is nonnegative. As $\alpha\to\infty$, $\widehat\Lambda_0(\alpha)\ge h_\alpha(0)=b(\alpha)-\tfrac12\sigma_3^2\to+\infty$. A continuous function on $(0,\infty)$ that is coercive at both ends attains its infimum.
\end{proof}

\begin{theorem}[Almost sure exponential extinction]\label{Syz}
Assume (H1)--(H3) and let $u$ solve \eqref{SM} with $u_0\in\mathbb{R}^3_+$. Then, almost surely,
\begin{equation}\label{eq:main_extinction}
\limsup_{t\to\infty}\ \frac1t\,\ln\big(y(t)+z(t)\big)\ \le\ \Lsharp:=\inf_{\alpha>0}\max\big\{\Lambda_\alpha,\ \Lambda_\alpha^{\flat}\big\}.
\end{equation}
If some minimizer $\alpha^\ast$ of \eqref{eq:lambda_ext_def} satisfies $(\mathrm{C}_{\alpha^\ast})$, then $\Lsharp=\Lext$ and the bound holds with $\Lext$ in place of $\Lsharp$. In particular, if $\Lsharp<0$ then $y(t)\to0$ and $z(t)\to0$ exponentially fast almost surely.
\end{theorem}

\begin{proof}
Fix $\alpha>0$; the final bound follows by taking the infimum over a countable dense subset of $(0,\infty)$ and collecting the exceptional null sets, the map $\alpha\mapsto\max\{\Lambda_\alpha,\Lambda_\alpha^\flat\}$ being continuous.

\textbf{Step 1 (It\^o's formula).} Let $S(t)=\alpha y(t)+z(t)$, which is strictly positive for all $t$ almost surely by Theorem~\ref{Existence}. Then
\begin{align}\label{eq:dV_ext}
d\ln S &= \frac{1}{S}\Big(\alpha\betaeta xz-\alpha(\mu_2+q)y+\peps y-\mu_3z\Big)dt-\frac{\alpha^2\sigma_2^2y^2+\sigma_3^2z^2}{2S^2}\,dt \nonumber\\
&\quad+\frac{\alpha\sigma_2y}{S}\,dW_2+\frac{\sigma_3z}{S}\,dW_3 .
\end{align}

\textbf{Step 2 (drift in barycentric form).} Put $w(t)=\alpha y(t)/S(t)\in(0,1)$, so that $1-w(t)=z(t)/S(t)$. Dividing each drift term by $S$ and collecting,
\begin{equation}\label{eq:LV_w}
\mathcal{L}\ln S = g_\alpha\big(w(t),x(t)\big),\qquad
g_\alpha(w,\chi):=w\,a(\alpha)+(1-w)\big(\alpha\betaeta\chi-\mu_3\big)-\tfrac12\sigma_2^2w^2-\tfrac12\sigma_3^2(1-w)^2 .
\end{equation}
Note that $g_\alpha(w,\bar x)=h_\alpha(w)$ and that, for each fixed $w\in[0,1]$, the map $\chi\mapsto g_\alpha(w,\chi)$ is affine and nondecreasing with slope $\alpha\betaeta(1-w)\ge0$.

\textbf{Step 3 (envelope in $x$).} Let $\Phi_\alpha(\chi):=\sup_{w\in[0,1]}g_\alpha(w,\chi)$. As a supremum of nondecreasing affine functions of $\chi$, $\Phi_\alpha$ is nondecreasing and convex; this convexity is precisely the obstruction to replacing $x$ by $\bar x$. Enlarging the range of $w$ to $\mathbb{R}$ and completing the square gives the explicit quadratic majorant
\begin{equation}\label{eq:Phihat}
\Phi_\alpha(\chi)\le\widehat\Phi_\alpha(\chi):=\alpha\betaeta\chi-\mu_3-\tfrac12\sigma_3^2+\frac{\big(a(\alpha)-\alpha\betaeta\chi+\mu_3+\sigma_3^2\big)^2}{2\varsigma^2}=A_2\chi^2+A_1\chi+A_0 ,
\end{equation}
with $A_2=A_2(\alpha)>0$ as in \eqref{eq:lambda_ext_def}, $A_0=\widehat\Phi_\alpha(0)$ as in \eqref{eq:Lflat}, and, since $\widehat\Phi_\alpha'(\bar x)=\alpha\betaeta\big(1-\widehat w(\alpha)\big)$,
\[
A_1=\alpha\betaeta\big(1-\widehat w(\alpha)\big)-2A_2\bar x .
\]
By construction $\widehat\Phi_\alpha(\bar x)=\widehat\Lambda_0(\alpha)$, and $\mathcal{L}\ln S(t)\le\widehat\Phi_\alpha(x(t))$ for every $t$.

\textbf{Step 4 (martingale term).} The local martingale $M(t)=\int_0^t\frac{\alpha\sigma_2y}{S}dW_2+\int_0^t\frac{\sigma_3z}{S}dW_3$ satisfies $d\langle M\rangle_t\le\varsigma^2dt$, since $\alpha y/S\le1$ and $z/S\le1$. Hence $\int_0^\infty(1+s)^{-2}d\langle M\rangle_s<\infty$ and \eqref{eq:SLLN} gives $M(t)/t\to0$ almost surely.

\textbf{Step 5 (linear  program over the occupation constraints).} Integrating \eqref{eq:dV_ext} and dividing by $t$,
\begin{equation}\label{eq:preLP}
\frac{\ln S(t)-\ln S(0)}{t}\ \le\ A_2\avg{x^2}+A_1\avg{x}+A_0+\frac{M(t)}{t}.
\end{equation}
Fix $\varepsilon>0$. Since $2\Lambda/(2\mu_1-\sigma_1^2)=\bar x(1+\rho_1)$, Lemma~\ref{lem:occupation} provides an almost surely finite $T_\varepsilon$ such that, for all $t\ge T_\varepsilon$,
\[
\avg{x^2}\le \bar x(1+\rho_1)\big(\avg{x}+\avg{y}\big)+\varepsilon,\qquad \mu_1\avg{x}+\mu_2\avg{y}\le\Lambda+\varepsilon .
\]
Because $A_2>0$, the first inequality gives
\[
A_2\avg{x^2}+A_1\avg{x}\ \le\ c_1\avg{x}+c_2\avg{y}+A_2\varepsilon,\qquad c_1:=A_2\bar x(1+\rho_1)+A_1,\quad c_2:=A_2\bar x(1+\rho_1).
\]
The pair $(\avg{x},\avg{y})$ lies in the simplex $\Delta_\varepsilon=\{(\xi,\upsilon):\xi\ge0,\ \upsilon\ge0,\ \mu_1\xi+\mu_2\upsilon\le\Lambda+\varepsilon\}$, whose vertices are $(0,0)$, $\big(\tfrac{\Lambda+\varepsilon}{\mu_1},0\big)$ and $\big(0,\tfrac{\Lambda+\varepsilon}{\mu_2}\big)$. A linear functional on a polytope attains its maximum at a vertex, so
\[
c_1\avg{x}+c_2\avg{y}\ \le\ \Big(1+\tfrac{\varepsilon}{\Lambda}\Big)\max\Big\{0,\ c_1\bar x,\ c_2\tfrac{\Lambda}{\mu_2}\Big\},
\]
the right-hand maximum being nonnegative. Letting $t\to\infty$ in \eqref{eq:preLP}, using Step 4, and then letting $\varepsilon\downarrow0$,
\[
\limsup_{t\to\infty}\frac{\ln S(t)}{t}\ \le\ A_0+\max\Big\{0,\ c_1\bar x,\ c_2\frac{\Lambda}{\mu_2}\Big\}\qquad\text{almost surely.}
\]
Using $\widehat\Lambda_0(\alpha)=A_2\bar x^2+A_1\bar x+A_0$ and $\Lambda=\mu_1\bar x$, the three candidate values are
\[
A_0=\widehat\Phi_\alpha(0),\qquad c_1\bar x+A_0=\widehat\Lambda_0(\alpha)+A_2\bar x^2\rho_1=\Lambda_\alpha,\qquad c_2\frac{\Lambda}{\mu_2}+A_0=\Lambda_\alpha^{\flat},
\]
and $\widehat\Phi_\alpha(0)\le\Lambda_\alpha^\flat$ because the term added in \eqref{eq:Lflat} is nonnegative. Since $\min\{\alpha,1\}(y+z)\le S\le\max\{\alpha,1\}(y+z)$, the same bound holds with $S$ replaced by $y+z$, and taking the infimum over $\alpha$ gives \eqref{eq:main_extinction}.

\textbf{Step 6 (when the penalty branch is inactive).} Substituting the expression for $A_1$ into $\Lambda_\alpha-\Lambda_\alpha^\flat=A_2\bar x^2(1+\rho_1)\big(1-\tfrac{\mu_1}{\mu_2}\big)+A_1\bar x$ gives
\begin{equation}\label{eq:Cgap}
\Lambda_\alpha-\Lambda_\alpha^{\flat}=\bar x\Big[\alpha\betaeta\big(1-\widehat w(\alpha)\big)-A_2(\alpha)\bar x\,\gamma\Big],
\end{equation}
so $\Lambda_\alpha^{\flat}\le\Lambda_\alpha$ precisely when $(\mathrm{C}_\alpha)$ holds, after division by $\alpha\betaeta>0$. If $(\mathrm{C}_{\alpha^\ast})$ holds at a minimizer $\alpha^\ast$, which exists by Lemma~\ref{lem:attain}, then $\Lsharp\le\max\{\Lambda_{\alpha^\ast},\Lambda_{\alpha^\ast}^\flat\}=\Lambda_{\alpha^\ast}=\Lext$, while $\Lsharp\ge\inf_\alpha\Lambda_\alpha=\Lext$; hence $\Lsharp=\Lext$.
\end{proof}

\subsection{The Deterministic Limit}

Theorem~\ref{Syz} gives a sufficient condition for clearance. The following definition fixes the precise sense in which it is sharp; no claim of necessity is made at positive noise, and no claim is made along every path to zero noise.

\begin{definition}[Sharpness of the reduced bound under the deterministic-limit scaling]
\label{def:sharp}
Call a sequence of noise intensities a deterministic-limit scaling if
$\sigma_2\to0$, $\sigma_3\to0$, and $\sigma_1^2/\varsigma^2\to0$; the case
$\sigma_1=0$ with $\sigma_2,\sigma_3\to0$ is included. We call the reduced
optimized bound $\Lext$ sharp under this scaling if
\[
\Lext\longrightarrow\lambda_{\max}(J)
\]
along every such sequence. Consequently, $\Lext<0$ reduces to
$\lambda_{\max}(J)<0$, equivalently $\Rz<1$, along this scaling. Nothing is
asserted along sequences for which $\sigma_1^2/\varsigma^2$ does not vanish.
\end{definition}

\begin{proposition}[Two-sided deterministic-limit estimate]\label{prop:detlimit}
Assume (H1)--(H3), and let $\alpha^\circ$ be the balancing weight of Lemma~\ref{lem:CW}. Then
\begin{equation}\label{eq:sandwich}
\lambda_{\max}(J)-\tfrac12\max\{\sigma_2^2,\sigma_3^2\}\ \le\ \Lext\ \le\ \lambda_{\max}(J)-\frac{\sigma_2^2\sigma_3^2}{2\varsigma^2}+A_2(\alpha^\circ)\,\bar x^{\,2}\rho_1 .
\end{equation}
In particular $\Lext\to\lambda_{\max}(J)$ along any sequence of noise intensities with $\sigma_2,\sigma_3\to0$ and $\sigma_1^2/\varsigma^2\to0$, and in particular whenever $\sigma_1=0$, so $\Lext$ is sharp under the deterministic-limit scaling of Definition~\ref{def:sharp}; if $\sigma_1=0$ then $\Lext<\lambda_{\max}(J)$ strictly whenever $\sigma_2\sigma_3>0$. The restriction $\sigma_1^2/\varsigma^2\to0$ cannot be dropped: at fixed $\sigma_1>0$ the penalty term forces $\Lext\to+\infty$ as $\varsigma\to0$, so the convention $\Lext:=\lambda_{\max}(J)$ at $\varsigma^2=0$ is a genuine discontinuity and the criterion is uninformative when the infected compartments are nearly noiseless but the hepatocyte pool is not. Moreover $\Lext$ is nonincreasing in each of $\sigma_2^2$ and $\sigma_3^2$ and nondecreasing in $\sigma_1^2$.
\end{proposition}

\begin{proof}
For the lower bound, the penalty in \eqref{eq:lambda_ext_def} is nonnegative, and
$\widehat\Lambda_0(\alpha)=\sup_{w\in\mathbb{R}}h_\alpha(w)$ dominates both $h_\alpha(1)=a(\alpha)-\tfrac12\sigma_2^2$ and $h_\alpha(0)=b(\alpha)-\tfrac12\sigma_3^2$, so
\[
\Lambda_\alpha\ \ge\ \max\{a(\alpha),b(\alpha)\}-\tfrac12\max\{\sigma_2^2,\sigma_3^2\},
\]
and taking the infimum over $\alpha$ and applying \eqref{eq:CW} gives the left inequality.

For the upper bound, evaluate at $\alpha^\circ$. There $a=b=\lambda_{\max}(J)$, so $\widehat w(\alpha^\circ)=\sigma_3^2/\varsigma^2$ and
\[
\widehat\Lambda_0(\alpha^\circ)=\lambda_{\max}(J)-\tfrac12\sigma_3^2+\frac{\sigma_3^4}{2\varsigma^2}=\lambda_{\max}(J)-\frac{\sigma_3^2(\varsigma^2-\sigma_3^2)}{2\varsigma^2}=\lambda_{\max}(J)-\frac{\sigma_2^2\sigma_3^2}{2\varsigma^2},
\]
and $\Lext\le\Lambda_{\alpha^\circ}=\widehat\Lambda_0(\alpha^\circ)+A_2(\alpha^\circ)\bar x^2\rho_1$. Both bounding expressions converge to $\lambda_{\max}(J)$ under the stated conditions, since $A_2(\alpha^\circ)\bar x^2\rho_1=\alpha^{\circ2}\betaeta^2\bar x^2\rho_1/(2\varsigma^2)$ and $\rho_1\le\sigma_1^2/\mu_1$ by \eqref{eq:rho1} and (H2), so the penalty is $O(\sigma_1^2/\varsigma^2)$; conversely it diverges as $\varsigma\to0$ at fixed $\sigma_1>0$, because $\alpha^\circ$ does not depend on the noise.

Monotonicity: for fixed $\alpha$ and each $w$, $h_\alpha(w)$ is nonincreasing in $\sigma_2^2$ and in $\sigma_3^2$, hence so is its supremum $\widehat\Lambda_0(\alpha)$; and $A_2(\alpha)=\alpha^2\betaeta^2/(2\varsigma^2)$ is decreasing in $\varsigma^2$. Thus each $\Lambda_\alpha$, and therefore their infimum, is nonincreasing in $\sigma_2^2$ and $\sigma_3^2$. Dependence on $\sigma_1$ enters only through $\rho_1$, which is increasing in $\sigma_1^2$.
\end{proof}

\begin{remark}[Mechanism, and the cost of hepatocyte fluctuations]\label{rem:extinction}
\emph{(a) Why the weight matters.} The unweighted choice $\alpha=1$, corresponding to the functional $\ln(y+z)$ used in earlier treatments, yields in the noise-free limit the criterion $\max\{a(1),b(1)\}<0$, that is $\peps<\mu_2+q$ and $\betaeta\bar x<\mu_3$ simultaneously. By Lemma~\ref{lem:CW} this is strictly stronger than $\Rz<1$ unless $\alpha^\circ=1$, and it fails by a wide margin at realistic parameters: in scenario S2 of Section~\ref{sec:numerics} it returns $\Lambda_1=+2.63$ day$^{-1}$ against $\Lext=-0.050$ day$^{-1}$.

\emph{(b) Role of $\sigma_2$ and $\sigma_3$.} The It\^o corrections enter through $-\tfrac12\sigma_2^2w^2-\tfrac12\sigma_3^2(1-w)^2$, so the barycentric coordinate $w$, the share of the infected load carried by hepatocytes rather than virions, determines how much of the available noise is felt, and the supremum over $w$ selects the least favorable allocation. The gain at the balancing weight, $\sigma_2^2\sigma_3^2/(2\varsigma^2)$ in \eqref{eq:sandwich}, is half the harmonic mean of $\sigma_2^2$ and $\sigma_3^2$: it vanishes if either compartment is noiseless, so within the model clearance is assisted only when both infected compartments fluctuate. The lower bound in \eqref{eq:sandwich} shows in addition that the total benefit cannot exceed $\tfrac12\max\{\sigma_2^2,\sigma_3^2\}$, which is why noise-induced clearance is a near-threshold phenomenon at biologically plausible intensities.

\emph{(c) Cost of hepatocyte fluctuations.} The penalty $A_2(\alpha)\bar x^2\rho_1$ is strictly positive when $\sigma_1>0$ and grows with the relative variance $\rho_1$ of the invariant law $\pi_1$. Fluctuations in the target-cell pool are thus an obstacle to clearance: because the escape rate is convex in $x$, excursions of $x$ above $\bar x$ contribute more than excursions below $\bar x$ subtract. Replacing $x(t)$ by $\bar x$ discards exactly this term. At the parameter values of Section~\ref{sec:numerics} the penalty is $4.35\times10^{-4}$ day$^{-1}$, numerically small but nonzero, and it grows quadratically in $\alpha\betaeta\bar x$.

\emph{(d) Status of the second branch.} The quantity $\Lambda_\alpha^{\flat}$ is the value taken by the same bound in the configuration where the hepatocyte pool is exhausted, $\avg{x}\to0$ and $\avg{y}\to\Lambda/\mu_2$. Condition $(\mathrm{C}_\alpha)$ states that this configuration is not the binding one; it is a restriction on the proof rather than on the dynamics, since the linear  program of Step 5 uses only two occupation constraints and is therefore blind to the fact that target-cell exhaustion is incompatible with clearance. Table~\ref{T3} reports $\Lambda_{\alpha^\ast}^\flat$ and the margin in $(\mathrm{C}_{\alpha^\ast})$; in the clearance regime studied here, the margin is positive and $\Lsharp=\Lext$.
\end{remark}

\subsection{Asymptotics of the Uninfected Compartment}

\begin{theorem}[Exponential convergence of $x$ to the reference process]\label{x-x1}
Assume the hypotheses of Theorem~\ref{Syz} together with $\Lsharp<0$, and let $x_1$ solve \eqref{x1-eq} with $x_1(0)=x_0=x(0)$. Then there exist a deterministic $\gamma_0>0$ and an almost surely finite random variable $C_\omega$ such that
\[
|x(t)-x_1(t)|\ \le\ C_\omega e^{-\gamma_0 t}\quad\text{for all }t\ge0;\qquad\text{in particular } \lim_{t\to\infty}|x(t)-x_1(t)|=0\ \text{ almost surely.}
\]
\end{theorem}

\begin{proof}
Let $e=x-x_1$. Subtracting \eqref{x1-eq} from the first equation of \eqref{SM},
\begin{equation}\label{eq:error_eq}
de(t)=\big[-\mu_1e(t)+F(t)\big]dt+\sigma_1e(t)\,dW_1(t),\qquad e(0)=0,\qquad F:=qy-\betaeta xz .
\end{equation}
The diffusion coefficient of \eqref{eq:error_eq} is exactly proportional to $e$, and the forcing $F$ has finite variation on compacts, so the variation of constants formula established in Lemma~\ref{lem4}(i) applies without an It\^o correction:
\begin{equation}\label{eq:varconst}
e(t)=\Phi(t)\int_0^t\Phi(s)^{-1}F(s)\,ds,\qquad \Phi(t)=\exp\Big(-\big(\mu_1+\tfrac12\sigma_1^2\big)t+\sigma_1W_1(t)\Big).
\end{equation}
Fix $\lambda\in(\Lsharp,0)$. By Theorem~\ref{Syz} and the continuity of the paths, there is an almost surely finite $C_1$ with $y(t)+z(t)\le C_1e^{\lambda t}$ for all $t\ge0$. By Lemma~\ref{lem:dissip}, for every $\varepsilon>0$ there is an almost surely finite $C_2$ with $x(t)\le W(u(t))\le C_2e^{\varepsilon t}$ for all $t\ge0$. Hence
\[
|F(t)|\le qC_1e^{\lambda t}+\betaeta C_1C_2e^{(\lambda+\varepsilon)t}\le C_3e^{(\lambda+\varepsilon)t},\qquad t\ge0 .
\]
Since $W_1(t)/t\to0$ almost surely, for every $\varepsilon>0$ the quantities $\sup_{t\ge0}\big(\sigma_1W_1(t)-\sigma_1\varepsilon t\big)$ and $\sup_{s\ge0}\big(-\sigma_1W_1(s)-\sigma_1\varepsilon s\big)$ are almost surely finite; writing $\alpha_0:=\mu_1+\tfrac12\sigma_1^2>0$, there is therefore an almost surely finite $C_4$ with
\[
\Phi(t)\le C_4e^{-(\alpha_0-\sigma_1\varepsilon)t},\qquad \Phi(s)^{-1}\le C_4e^{(\alpha_0+\sigma_1\varepsilon)s},\qquad t,s\ge0 .
\]
Choose $\varepsilon>0$ so small that $\alpha_0-\sigma_1\varepsilon>0$ and $\varepsilon(1+2\sigma_1)<|\lambda|$, and set $\theta_0:=-(\lambda+\varepsilon)>0$, so that $\theta_0-2\sigma_1\varepsilon>0$. Using $\int_0^te^{\kappa s}ds\le(1+t)e^{\kappa^+t}$ for every $\kappa\in\mathbb{R}$, with $\kappa=\alpha_0+\sigma_1\varepsilon-\theta_0$,
\[
|e(t)|\ \le\ C_3C_4^2\,(1+t)\,e^{-(\alpha_0-\sigma_1\varepsilon)t}e^{(\alpha_0+\sigma_1\varepsilon-\theta_0)^+t}\ \le\ C_5(1+t)e^{-2\gamma_0t},
\]
where $2\gamma_0:=\min\{\alpha_0-\sigma_1\varepsilon,\ \theta_0-2\sigma_1\varepsilon\}>0$: indeed if $\kappa\le0$ the exponent is $-(\alpha_0-\sigma_1\varepsilon)$, and if $\kappa>0$ it is $2\sigma_1\varepsilon-\theta_0$. Absorbing the factor $1+t$ into the exponential gives $|e(t)|\le C_\omega e^{-\gamma_0t}$ with $C_\omega=C_5\sup_{t\ge0}(1+t)e^{-\gamma_0t}<\infty$.
\end{proof}

\begin{corollary}[Long-run behaviour under clearance]\label{cor:attained}
Under the hypotheses of Theorem~\ref{x-x1}, almost surely:
\begin{enumerate}
\item[(i)] for every bounded uniformly continuous $g$ on $(0,\infty)$, $\avg{g(x)}\to\int g\,d\pi_1$;
\item[(ii)] $\avg{x^2}\to\bar x^{\,2}(1+\rho_1)$, so the occupation bound \eqref{eq:secondmomentbound} is attained and cannot be improved;
\item[(iii)] $u(t)$ converges in distribution to $\pi_1\otimes\delta_0\otimes\delta_0$ on $\overline{\mathbb{R}^3_+}$.
\end{enumerate}
\end{corollary}

\begin{proof}
(i) By Theorem~\ref{x-x1}, $x(s)-x_1(s)\to0$ almost surely; uniform continuity of $g$ gives $g(x(s))-g(x_1(s))\to0$, and boundedness allows the Ces\`aro limit to be taken, so $\avg{g(x)}-\avg{g(x_1)}\to0$. Lemma~\ref{lem4}(iii) supplies $\avg{g(x_1)}\to\int g\,d\pi_1$, $g$ being bounded and hence $\pi_1$-integrable.

(ii) Write $|x^2-x_1^2|=|x-x_1|(x+x_1)$. Theorem~\ref{x-x1} gives $|x-x_1|\le C_\omega e^{-\gamma_0s}$, and since $x_1\le x+C_\omega$ and $x\le C_2e^{\varepsilon s}$ with $\varepsilon<\gamma_0$ by Lemma~\ref{lem:dissip}, the product is almost surely integrable on $[0,\infty)$, so $\avg{x^2}-\avg{x_1^2}\to0$. Under (H2) one has $\sigma_1^2<2\mu_1$, so $\chi\mapsto\chi^2$ is $\pi_1$-integrable and Lemma~\ref{lem4}(iii) gives $\avg{x_1^2}\to\mathbb{E}_{\pi_1}[\chi^2]=\bar x^2(1+\rho_1)$ by \eqref{eq:secondmoment}. Since \eqref{eq:secondmomentbound} holds with the same constant, no smaller bound is valid for all parameter sets.

(iii) $x_1(t)\Rightarrow\pi_1$ by Lemma~\ref{lem4}(iii), and $x(t)-x_1(t)\to0$ almost surely, hence in probability, so $x(t)\Rightarrow\pi_1$ by Slutsky's theorem; $y(t),z(t)\to0$ almost surely by Theorem~\ref{Syz}. Joint convergence follows because the limits of $y$ and $z$ are deterministic.
\end{proof}

\begin{remark}[Status of the limiting measure]\label{rem:boundary}
The measure $\pi_1\otimes\delta_0\otimes\delta_0$ appearing in Corollary~\ref{cor:attained}(iii) is \emph{not} an invariant probability measure of \eqref{SM} on the open orthant $\mathbb{R}^3_+$, on which the process lives for all finite times by Theorem~\ref{Existence}. It is the invariant law of the boundary dynamics of the natural extension of \eqref{SM} to the closed orthant $\overline{\mathbb{R}^3_+}$, restricted to the invariant face $\{(\chi,0,0):\chi>0\}$ on which \eqref{SM} reduces to \eqref{x1-eq}. What Corollary~\ref{cor:attained} asserts is weak convergence of the one-dimensional laws together with almost sure convergence of the empirical averages, and nothing stronger.
\end{remark}

\section{Ergodic Stationary Distribution under Persistence}\label{sec:ergodic}

We now treat the complementary situation, in which the noise is not strong enough to eliminate the infection, and the process settles into a nondegenerate stochastic endemic state. Since Theorem~\ref{Syz} forces two coordinates to zero when $\Lsharp<0$, ergodicity on $\mathbb{R}^3_+$ can hold only under hypotheses incompatible with that condition; Proposition~\ref{prop:dichotomy} confirms that this is so.

\begin{lemma}[Positive recurrence and ergodicity {\citep[Thm.~4.1]{Ref29}}, {\citep[Lem.~2.3]{Ref33}}]\label{lem:MT}
Let $U(t)$ be a regular time-homogeneous Markov diffusion on an open set $\mathcal{O}\subseteq\mathbb{R}^d$ with generator $\mathcal{L}$ and diffusion matrix $A(u)=B(u)B(u)^\top$. Suppose there exist a nonempty bounded open set $D$ with $\overline{D}\subset\mathcal{O}$ and a function $V\in C^2(\mathcal{O})$ with $V\ge1$ and $V(u)\to\infty$ as $u\to\partial\mathcal{O}$ or $|u|\to\infty$, such that
\begin{align}
&\text{\emph{(A1)}}\quad \exists\,\underline{\lambda}>0:\ \xi^{\top}A(u)\xi\ge \underline{\lambda}|\xi|^2\quad \forall u\in D,\ \forall\xi\in\mathbb{R}^d; \label{A1}\\
&\text{\emph{(A2)}}\quad \exists\,\vartheta>0:\ \mathcal{L}V(u)\le -\vartheta \quad \forall u\in \mathcal{O}\setminus D. \label{A2}
\end{align}
Then $U$ is positive recurrent, admits a unique invariant probability measure $\pi$ on $\mathcal{O}$, and is ergodic, in the sense that $\frac1T\int_0^Tg(U(t))\,dt\to\int g\,d\pi$ almost surely for every $\pi$-integrable $g$.
\end{lemma}

Define
\begin{equation}\label{eq:mi}
m_1=\mu_1+\tfrac12\sigma_1^2,\qquad m_2=\mu_2+q+\tfrac12\sigma_2^2,\qquad m_3=\mu_3+\tfrac12\sigma_3^2,
\end{equation}
and
\begin{equation}\label{eq:Rs}
\Rs\ :=\ \frac{(1-\eta)(1-\epsilon)\beta p\Lambda}{m_1m_2m_3}\ =\ \frac{\betaeta\peps\Lambda}{\big(\mu_1+\frac{\sigma_1^2}{2}\big)\big(\mu_2+q+\frac{\sigma_2^2}{2}\big)\big(\mu_3+\frac{\sigma_3^2}{2}\big)} .
\end{equation}
Each $m_i$ is the removal rate of compartment $i$ augmented by the It\^o correction of the corresponding logarithmic Lyapunov function, and $\Rs$ is precisely the ratio appearing in the arithmetic--geometric mean balance of the proof below. When $\sigma_1=\sigma_2=\sigma_3=0$ it reduces to $\Rz$, and it is strictly decreasing in each $\sigma_i$. The correction to $\mu_1$ is not optional: omitting it produces a threshold that the Foster--Lyapunov argument does not deliver.

\begin{theorem}[Unique ergodic stationary distribution]\label{thm:ergodic}
Assume $\sigma_1,\sigma_2,\sigma_3>0$ and $\Rs>1$. Then the Markov process $U(t)=(x(t),y(t),z(t))$ generated by \eqref{SM} is positive recurrent on $\mathbb{R}^3_+$ and admits a unique invariant probability measure $\pi$ on $\mathbb{R}^3_+$; it is ergodic in the sense of Lemma~\ref{lem:MT}.
\end{theorem}

\begin{proof}
\emph{Step 1 (ellipticity on compacta).} The diffusion matrix is $A(u)=\mathrm{diag}(\sigma_1^2x^2,\sigma_2^2y^2,\sigma_3^2z^2)$. On any set $D$ whose closure is a compact subset of $\mathbb{R}^3_+$ all coordinates are bounded below by some $\varrho_0>0$, so $\xi^\top A(u)\xi\ge\varrho_0^2\min_i\sigma_i^2\,|\xi|^2$ and \eqref{A1} holds. Regularity of the diffusion follows from Theorem~\ref{Existence}.

\emph{Step 2 (Lyapunov function).} Put $c_i=1/m_i$ for $i=1,2,3$ and
\[
V_1(u)=-c_1\ln x-c_2\ln y-c_3\ln z,\qquad V_\theta(u)=\tfrac{1}{\theta+1}W(u)^{\theta+1},
\]
with $W$, $\delta_1$, $\delta_2$, $c_z$ and $\theta$ as in Lemma~\ref{lem:dissip}, and set $V=MV_1+V_\theta+C$ with $M>0$ and $C$ chosen below. Writing $m:=\min\{1,c_z\}$ we have $\max_iu_i\le W/m$, so
\[
MV_1(u)+V_\theta(u)\ \ge\ -M\Big(\textstyle\sum_ic_i\Big)\ln\frac{W}{m}+\frac{W^{\theta+1}}{\theta+1},
\]
whose right-hand side is a function of $W$ alone tending to $+\infty$ as $W\to0^+$ and as $W\to\infty$, hence bounded below; choosing $C$ large makes $V\ge1$. If $u\to\partial\mathbb{R}^3_+$ then some $u_i\to0^+$ while $W$ stays bounded away from $0$ along the same sequence, or else the displayed bound already diverges, so $V\to\infty$; and $V\to\infty$ as $|u|\to\infty$ because $W^{\theta+1}$ dominates $\ln W$.

\emph{Step 3 (generator of $V_1$).} Using $\partial_v(-\ln v)=-1/v$ and $\partial_v^2(-\ln v)=1/v^2$, each coordinate contributes its It\^o correction $\tfrac12\sigma_i^2$, so
\[
\mathcal{L}V_1(u)=-G(u)+c_1m_1+c_2m_2+c_3m_3+c_1\betaeta z-\frac{c_1qy}{x}\ \le\ -G(u)+3+c_1\betaeta z,
\]
where
$G(u)=\dfrac{c_1\Lambda}{x}+\dfrac{c_2\betaeta xz}{y}+\dfrac{c_3\peps y}{z}$
and $c_1m_1+c_2m_2+c_3m_3=3$. The three summands of $G$ have product
\[
\frac{c_1\Lambda}{x}\cdot\frac{c_2\betaeta xz}{y}\cdot\frac{c_3\peps y}{z}=c_1c_2c_3\Lambda\betaeta\peps=\frac{\Lambda\betaeta\peps}{m_1m_2m_3}=\Rs
\]
independently of $u$, every occurrence of $x$, $y$ and $z$ cancelling, so the arithmetic--geometric mean inequality gives the uniform bound $G(u)\ge3(\Rs)^{1/3}$. Fix $\epsilon_0\in(0,1)$ so small that $(1-\epsilon_0)(\Rs)^{1/3}>1$, which is possible because $\Rs>1$, and set $\lambda_1:=3\big[(1-\epsilon_0)(\Rs)^{1/3}-1\big]>0$. Splitting $G=\epsilon_0G+(1-\epsilon_0)G$ and applying the bound only to the second part,
\begin{equation}\label{eq:LV1}
\mathcal{L}V_1(u)\ \le\ -\epsilon_0G(u)-\lambda_1+c_1\betaeta z .
\end{equation}

\emph{Step 4 (choice of constants).} Let $C_3:=\sup_{r>0}\big(\Lambda r^\theta-\delta_2r^{\theta+1}\big)<\infty$. Choose, in this order:
\begin{enumerate}
\item[(i)] $M>0$ with $M\lambda_1\ge C_3+2$;
\item[(ii)] $\zeta:=\big(Mc_1\betaeta\big)^{-1}>0$, so that $Mc_1\betaeta z\le1$ whenever $z\le\zeta$;
\item[(iii)] $R_1>1+c_z\zeta$ so large that $\dfrac{Mc_1\betaeta}{c_z}r+\Lambda r^\theta-\delta_2r^{\theta+1}\le-1$ for all $r\ge R_1$, which is possible since $\theta+1>1$;
\item[(iv)] $\varrho\in(0,\tfrac12)$ so small that $M\epsilon_0c_1\Lambda/\varrho\ \ge\ Mc_1\betaeta R_1/c_z+\Lambda R_1^{\theta}+1$;
\item[(v)] $\varrho'\in(0,\varrho)$ so small that $M\epsilon_0c_2\betaeta\varrho\zeta/\varrho'\ \ge\ Mc_1\betaeta R_1/c_z+\Lambda R_1^{\theta}+1$.
\end{enumerate}
Each choice depends only on quantities already fixed, so no circularity arises. Define the bounded open set.
\[
D=\big\{u\in\mathbb{R}^3_+:\ x>\varrho,\ y>\varrho',\ z>\zeta,\ W(u)<R_1\big\},
\]
whose closure is a compact subset of $\mathbb{R}^3_+$: the constraint $W<R_1$ bounds all coordinates above, and the three strict inequalities bound them away from zero. It is nonempty, since $\varrho+\varrho'+c_z\zeta<1+c_z\zeta<R_1$.

\emph{Step 5 (drift condition outside $D$).} By \eqref{eq:LV1} and \eqref{eq:LVtheta},
\begin{equation}\label{eq:LVfull}
\mathcal{L}V(u)\ \le\ -M\epsilon_0G(u)-M\lambda_1+Mc_1\betaeta z+\Lambda W^\theta-\delta_2W^{\theta+1}.
\end{equation}
We verify $\mathcal{L}V\le-1$ on each part of $\mathbb{R}^3_+\setminus D$, using $G\ge0$ and $M\lambda_1>0$ throughout.
\begin{itemize}
\item $W\ge R_1$: since $z\le W/c_z$, item (iii) applied with $r=W$ gives $Mc_1\betaeta z+\Lambda W^\theta-\delta_2W^{\theta+1}\le-1$, so $\mathcal{L}V\le-M\lambda_1-1\le-1$.
\item $W<R_1$ and $z\le\zeta$: item (ii) gives $Mc_1\betaeta z\le1$, and $\Lambda W^\theta-\delta_2W^{\theta+1}\le C_3$, so item (i) gives $\mathcal{L}V\le-M\lambda_1+1+C_3\le-1$.
\item $W<R_1$, $z>\zeta$ and $x\le\varrho$: then $G(u)\ge c_1\Lambda/\varrho$, while $Mc_1\betaeta z\le Mc_1\betaeta R_1/c_z$ and $\Lambda W^\theta\le\Lambda R_1^\theta$; item (iv) gives $\mathcal{L}V\le-1$.
\item $W<R_1$, $z>\zeta$, $x>\varrho$ and $y\le\varrho'$: then $G(u)\ge c_2\betaeta xz/y\ge c_2\betaeta\varrho\zeta/\varrho'$, and item (v) gives $\mathcal{L}V\le-1$.
\end{itemize}
The four cases are stated as conditions on regions, not as limits, and together they exhaust $\mathbb{R}^3_+\setminus D$; in particular, states in which several coordinates are simultaneously small are covered by whichever of the last three cases applies first in the stated order. Hence \eqref{A2} holds with $\vartheta=1$.

\emph{Step 6 (conclusion).} Lemma~\ref{lem:MT} applies with $\mathcal{O}=\mathbb{R}^3_+$, giving positive recurrence, a unique invariant probability measure $\pi$ on $\mathbb{R}^3_+$, and ergodicity.
\end{proof}

\begin{proposition}[The two regimes are mutually exclusive]\label{prop:dichotomy}
Assume (H3). If $\Lext<0$ then $\Rs<1$. Consequently, the hypotheses of Theorems~\ref{Syz} and \ref{thm:ergodic} cannot hold simultaneously.
\end{proposition}

\begin{proof}
If $\Lext<0$ there is $\alpha>0$ with $\Lambda_\alpha<0$, hence $\widehat\Lambda_0(\alpha)<0$ because the penalty term in \eqref{eq:lambda_ext_def} is nonnegative. Since $\widehat\Lambda_0(\alpha)=\sup_{w\in\mathbb{R}}h_\alpha(w)\ge\max\{h_\alpha(0),h_\alpha(1)\}$, with $h_\alpha(1)=a(\alpha)-\tfrac12\sigma_2^2=\peps/\alpha-m_2$ and $h_\alpha(0)=b(\alpha)-\tfrac12\sigma_3^2=\alpha\betaeta\bar x-m_3$, we obtain
\[
\frac{\peps}{\alpha}<m_2,\qquad \alpha\betaeta\bar x<m_3 .
\]
All four quantities are positive, so multiplying eliminates $\alpha$ and yields $\peps\betaeta\bar x<m_2m_3$, that is $\betaeta\peps\Lambda<\mu_1m_2m_3$. Since $m_1>\mu_1$, this gives $\Rs=\betaeta\peps\Lambda/(m_1m_2m_3)<1$. The same conclusion follows a fortiori from $\Lsharp<0$, because $\Lsharp\ge\Lext$.
\end{proof}

\begin{remark}[The intermediate regime]\label{rem:separation}
Proposition~\ref{prop:dichotomy} shows that $\{\Lext<0\}\subset\{\Rs<1\}$, and the inclusion is strict: there are parameter sets with $\Rs<1$ and $\Lext\ge0$ for which neither theorem applies. In that band, the infection is not certified to clear at an exponential rate, and no nondegenerate invariant law is certified either; the present sufficient criteria leave this unresolved. Since Theorem~\ref{Existence} gives strict positivity of all three components at every finite time, the process is not absorbed at the boundary in finite time, so whatever occurs in this band is not absorption. Characterizing the regime may require quasi-stationary or related boundary-asymptotic methods; we do not establish that such a description applies to \eqref{SM}. Establishing which of the two outcomes occurs and locating the exact critical surface remain open. Figure~\ref{fig:regime} displays the three regions numerically.
\end{remark}

\section{Numerical Experiments}\label{sec:numerics}

\subsection{Scheme, Parameters, and Scenarios}

System \eqref{SM} has no closed-form solution, so we discretize it. Let $t_n=n\Delta t$, $n=0,\dots,N$, with $T=N\Delta t$, and let $\Delta W_{i,n}\sim\mathcal{N}(0,\Delta t)$ be independent increments. The Euler--Maruyama update is
\begin{align}
\tilde{x}_{n+1} &= x_n + \big(\Lambda-\mu_1x_n-\betaeta x_nz_n+qy_n\big)\Delta t+\sigma_1x_n\Delta W_{1,n}, \label{eq:EM_x}\\
\tilde{y}_{n+1} &= y_n + \big(\betaeta x_nz_n-(\mu_2+q)y_n\big)\Delta t+\sigma_2y_n\Delta W_{2,n}, \label{eq:EM_y}\\
\tilde{z}_{n+1} &= z_n + \big(\peps y_n-\mu_3z_n\big)\Delta t+\sigma_3z_n\Delta W_{3,n}, \label{eq:EM_z}
\end{align}
followed by the projection
\begin{equation}\label{eq:projection}
x_{n+1}=\max(\tilde{x}_{n+1},0),\qquad y_{n+1}=\max(\tilde{y}_{n+1},0),\qquad z_{n+1}=\max(\tilde{z}_{n+1},0).
\end{equation}
The projection is required because the exact solution is strictly positive by Theorem~\ref{Existence}. In contrast, the discretization is not: a sufficiently negative Gaussian increment can push a component below zero with positive probability at every step. It is the metric projection onto the closed positive cone; hence, it is the least intrusive repair available, and it does not move an admissible step. Two consequences should be recorded. First, the classical strong order $1/2$ theory \citep{kloeden2013numerical} covers neither the non-globally-Lipschitz drift of \eqref{SM} nor the non-smooth map \eqref{eq:projection}; we make no claim of a proved order for the projected scheme, and treat the order as a quantity to be measured. Methods that preserve positivity exactly, such as tamed or Lamperti-transformed schemes \citep{higham2001algorithmic}, would be preferable where the projection activates frequently, and the accompanying script records how often it does so. Second, once a step drives $y$ or $z$ to zero, the projected scheme is absorbed there, whereas the exact solution is not; all decay rates fitted below are therefore computed on the window in which the numerical trajectory is still strictly positive.

Baseline parameters are those of Table~\ref{T1}. They give $\bar x=7.14\times10^{8}$ cells and $\Rz=91.6$, so the untreated system lies far above the deterministic threshold, and this governs the interpretation of every simulation below: at $\Rz=91.6$ one has $\lambda_{\max}(J)=3.91$ day$^{-1}$, and solving $\Lext(\sigma)=0$ with $\sigma_2=\sigma_3=\sigma$ gives $\sigma\approx4.0$ day$^{-1/2}$. Noise of that magnitude is not biologically defensible, so noise-induced clearance is a near-threshold phenomenon here, as Proposition~\ref{prop:detlimit} predicts, and we study it under strong combination therapy, where it is relevant.

We therefore use the three scenarios of Table~\ref{T3}: an endemic scenario at baseline treatment, a near-threshold scenario in which moderate noise produces clearance, and the corresponding noise-free control. In the near-threshold scenario, the cellular noise $\sigma_2=0.35$ lies inside the bound $\sqrt{2\mu_2/3}=0.4$ imposed by (H2), while the hypotheses unconstrain the virion noise $\sigma_3=0.50$ and carry most of the clearance effect. The standing hypotheses hold in every scenario with $\varsigma^2>0$, since $\mu_1=0.014\le\mu_2=0.24$ and $\sigma_1=0.01<\sqrt{2\mu_1/3}=0.0966$; scenario S3 has $\varsigma^2=0$ and serves only as the deterministic comparison, for which $\Lext=\lambda_{\max}(J)$ by convention. Trajectory runs use $\Delta t=10^{-3}$ days and the long-run statistics $\Delta t=10^{-2}$ days, with initial state $(10^{7},10^{5},10^{4})$ throughout and one fixed random stream per trajectory, so every figure is reproducible bit for bit; the accompanying script regenerates all seven figures and prints the quantities of Table~\ref{T3}.

\begin{table}[H]
\centering
\caption{Simulation scenarios and the associated threshold quantities, computed from \eqref{R0}, \eqref{eq:lambda_ext_def}, \eqref{eq:Lflat} and \eqref{eq:Rs}, with $\sigma_1=0.01$ except in S3. The last two columns report the second branch of Theorem~\ref{Syz} at the optimal weight and the margin $1-\widehat w(\alpha^\ast)-\gamma\alpha^\ast\betaeta\bar x/(2\varsigma^2)$ in condition $(\mathrm{C}_{\alpha^\ast})$; a positive margin certifies $\Lsharp=\Lext$. All unlisted parameters are as in Table~\ref{T1}.}
\small
\begin{tabular}{lccccccccc}
\toprule
\textbf{Scenario} & $\eta$ & $\epsilon$ & $\sigma_2$ & $\sigma_3$ & $\Rz$ & $\Rs$ & $\Lext$ & $\Lambda^{\flat}_{\alpha^\ast}$ & $(\mathrm{C}_{\alpha^\ast})$ \\
\midrule
S1 endemic                 & $0.36$ & $0.24$ & $0.05$ & $0.05$ & $91.630$ & $90.706$ & $+10.082$ & $+1829$ & $-438$ \\
S2 noise-induced clearance & $0.94$ & $0.90$ & $0.35$ & $0.50$ & $1.130$  & $0.694$  & $-0.0504$ & $-0.1422$ & $+0.305$ \\
S3 noise-free control      & $0.94$ & $0.90$ & $0$    & $0$    & $1.130$  & $1.130$  & $+0.0242$ & --- & --- \\
\bottomrule
\end{tabular}
\label{T3}
\end{table}

Scenario S2 is the most informative case. Here, $\Rz=1.130>1$, so the deterministic model predicts persistence, as confirmed by the noise-free control S3, for which $\Lext=\lambda_{\max}(J)=+0.0242$ day$^{-1}$. Introducing stochastic fluctuations reduces the optimized bound to $\Lext=-0.0504$ day$^{-1}$. Since condition $(\mathrm{C}_{\alpha^\ast})$ is satisfied, $\Lsharp=\Lext$, and Theorem~\ref{Syz} therefore guarantees almost sure exponential clearance, with the asymptotic exponent of the infected population bounded above by $-0.0504$ day$^{-1}$. The optimal weight is $\alpha^\ast=0.292$, and the hepatocyte fluctuation penalty contributes only $4.35\times10^{-4}$ day$^{-1}$, less than one percent of $|\Lext|$. By contrast, the unweighted choice $\alpha=1$ gives $\Lambda_1=+2.63$ day$^{-1}$ and therefore fails to certify clearance, demonstrating the substantial improvement obtained by optimizing over $\alpha$. More specifically, the margin in $(\mathrm{C}_{\alpha^\ast})$ is $+0.305$, while $\Lambda^\flat_{\alpha^\ast}=-0.142<\Lext$, confirming that the general bound reduces to $\Lsharp=\Lext$ in this case. In the endemic scenario S1, the admissibility margin is strongly negative; however, since $\Lext=+10.1>0$, the extinction criterion is already inconclusive, so failure of the admissibility condition does not affect the classification of this scenario.
\subsection{Extinction, Persistence, and the Role of Noise}

\begin{figure}[!htbp]
\centering
\includegraphics[width=0.97\textwidth]{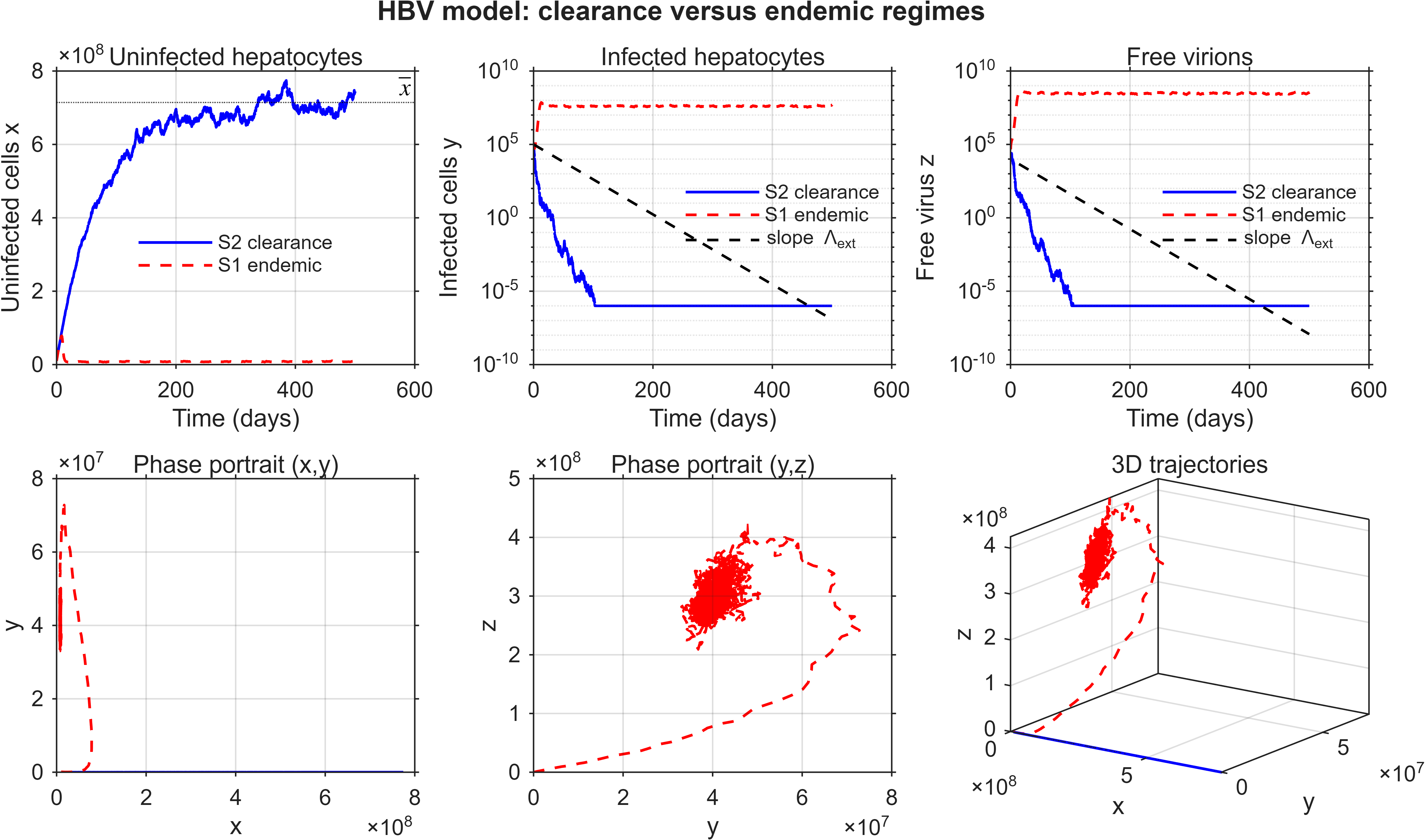}
\caption{\textbf{Clearance and endemic regimes under a common noise realization.} Top row: trajectories of $x(t)$, $y(t)$, and $z(t)$ for S2 (clearance, solid) and S1 (endemic, dashed). The reference line has a slope $\Lext=-0.0504$ day$^{-1}$; the flattening of the S2-infected compartments near $10^{-6}$ reflects numerical absorption in the projected scheme. Bottom row: corresponding phase-space projections and the three-dimensional trajectories. In S2, the infected compartments decay with an exponent below the analytical bound $\Lext$, while $x$ fluctuates about $\bar{x}=\Lambda/\mu_1$; in S1, the trajectory remains in the interior, consistent with the ergodic invariant measure of Theorem~\ref{thm:ergodic}.}
\label{fig:df-end}
\end{figure}

Figure~\ref{fig:df-end} contrasts the two regimes using the same Brownian realization, so the differences in the sample paths reflect the parameter regimes rather than different noise realizations. In scenario S2, condition $(\mathrm{C}_{\alpha^\ast})$ holds, so $\Lsharp=\Lext$ and Theorem~\ref{Syz} guarantees that the asymptotic decay rate of $\ln(y+z)$ does not exceed $\Lext$. The fitted rate is consistent with this one-sided bound and is more negative than $\Lext$. The gap is expected because several estimates in the proof are conservative: Step~3 replaces $\Phi_\alpha$ by the upper envelope $\widehat{\Phi}_\alpha$, and Step~5 relaxes the dynamics to an occupation-measure linear program. The realized infected-mass fraction $w(t)$ need not attain its worst-case allocation. Thus, $\Lext$ should be interpreted as a rigorous upper bound on the asymptotic exponent rather than an exact prediction of the observed decay rate.

\begin{figure}[!htbp]
\centering
\includegraphics[width=0.92\textwidth]{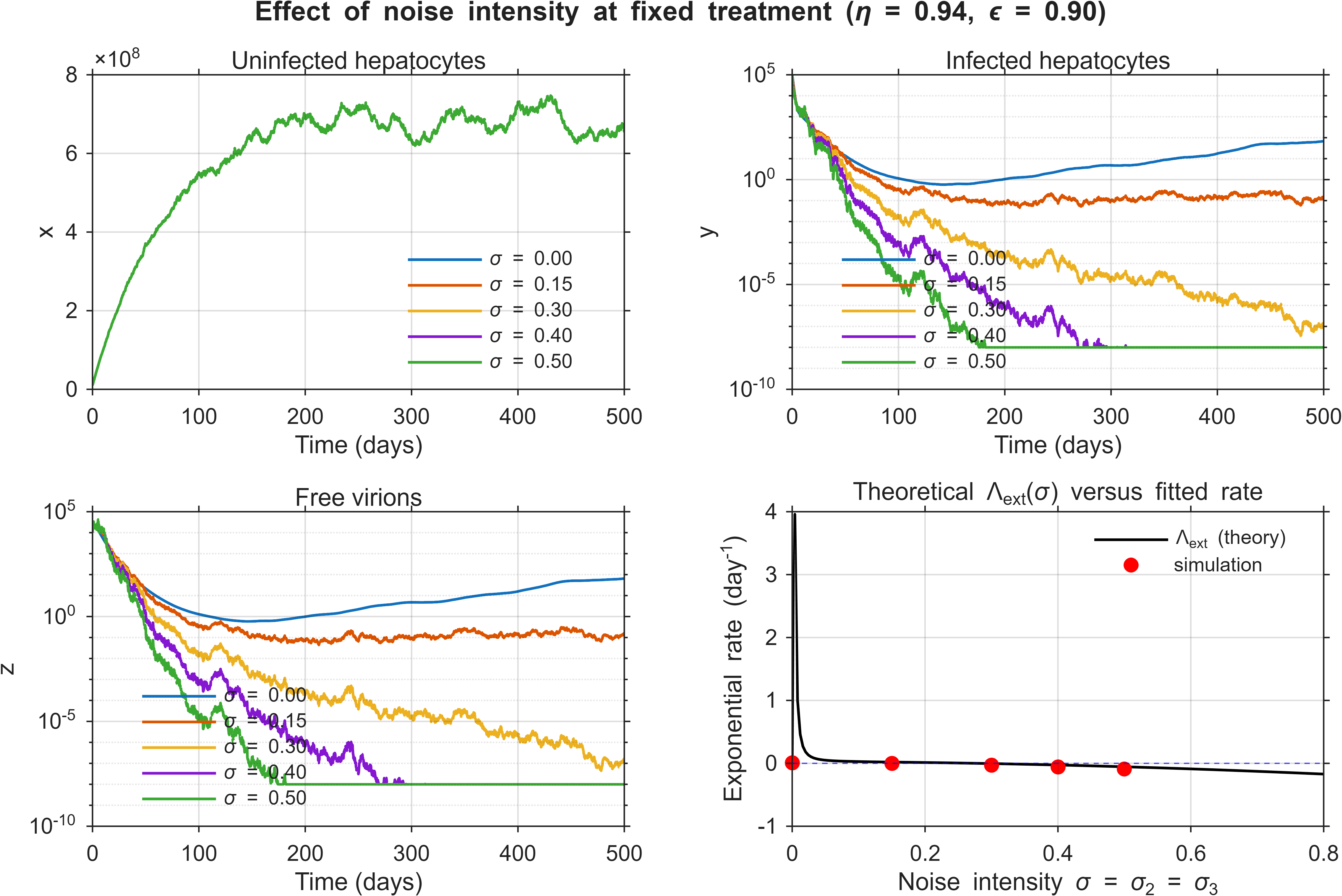}
\caption{\textbf{Effect of virion noise at fixed treatment}, with $\eta=0.94$, $\epsilon=0.90$, $\sigma_1=0.01$, $\sigma_2=0.35$, and $\sigma_3\in\{0,0.15,0.30,0.45,0.60\}$. The first three panels show $x(t)$, $y(t)$, and $z(t)$; flat tails indicate numerical absorption by the projected scheme, while the exact solution remains positive. The fourth compares fitted decay rates of $\ln(y+z)$ with the reduced optimized bound $\Lext(\sigma_3)$, which changes sign at $\sigma_3^\ast=0.269$. Where $(\mathrm{C}_{\alpha^\ast})$ holds, $\Lext<0$ certifies almost sure exponential clearance, and the fitted rates remain below the analytical bound, consistent with Theorem~\ref{Syz}.}
\label{fig:noise}
\end{figure}

Figure~\ref{fig:noise} isolates the effect of noise at fixed drift parameters. Along the sweep the analytical threshold takes the values $\Lext=+0.0223$, $+0.0153$, $-0.0052$, $-0.0374$ and $-0.0789$ day$^{-1}$, with $\Rs=1.055$, $1.008$, $0.889$, $0.742$ and $0.603$; the sign change occurs at $\sigma_3^\ast=0.269$. The monotone decrease is the content of the last statement of Proposition~\ref{prop:detlimit}. Its mechanism is the one identified in Remark~\ref{rem:extinction}(b): the It\^o corrections lower $h_\alpha$ uniformly in $w$, and once the maximum of $h_{\alpha^\ast}$ falls below zero, no allocation of infected mass between hepatocytes and virions can sustain growth. The transition is abrupt in the certified exponent, the quantity the theorem controls, but not in the sample paths, which cross over smoothly; and the fitted rate turns negative before $\Lext$ does, the interval between the two sign changes lying inside the intermediate band of Remark~\ref{rem:separation}.

\begin{figure}[!htbp]
\centering
\includegraphics[width=\textwidth]{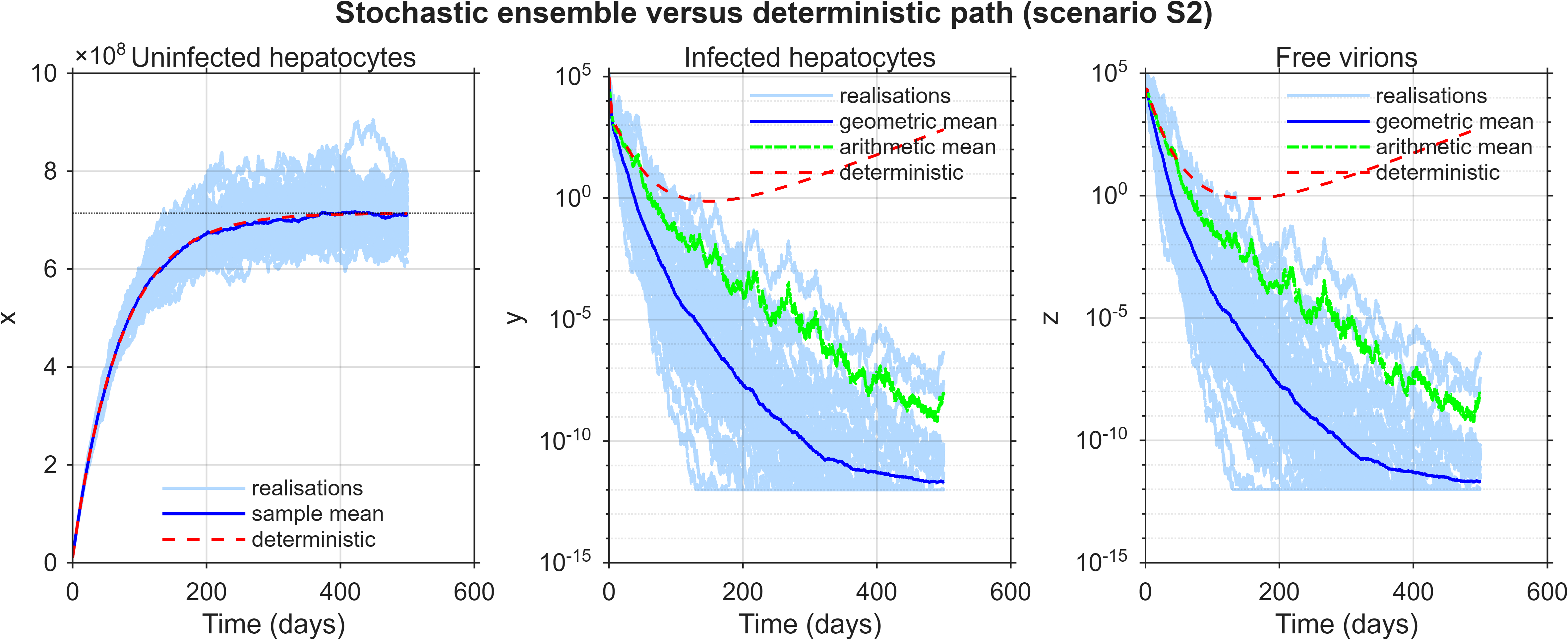}
\caption{\textbf{Ensemble against deterministic path.} Fifty independent realizations of scenario S2 (light curves) with, from left to right, $x(t)$, $y(t)$ and $z(t)$. The first panel shows the sample mean; the second and third show both the geometric mean, which tracks the typical path, and the arithmetic mean, which does not. The dashed curve is the deterministic solution obtained by setting $\sigma_i=0$, that is, scenario S3. The deterministic infected load grows while every stochastic realization decays, and the arithmetic mean decays strictly more slowly than the geometric mean, illustrating that noise-induced clearance is an almost-sure, not a mean-value, phenomenon.}
\label{fig:ensemble}
\end{figure}

Figure~\ref{fig:ensemble} makes a point that is easy to miss. Almost sure exponential extinction is a statement about $t^{-1}\ln(y+z)$, not about $\mathbb{E}[y+z]$, and the two can disagree: under multiplicative noise, the first moment may decay strictly more slowly than the typical path, or not at all, being dominated by rare large excursions. The separation between the geometric and arithmetic means is exactly this effect, and it is why moment-stability criteria are not equivalent to Theorem~\ref{Syz}.

\begin{figure}[!htbp]
\centering
\includegraphics[width=\textwidth]{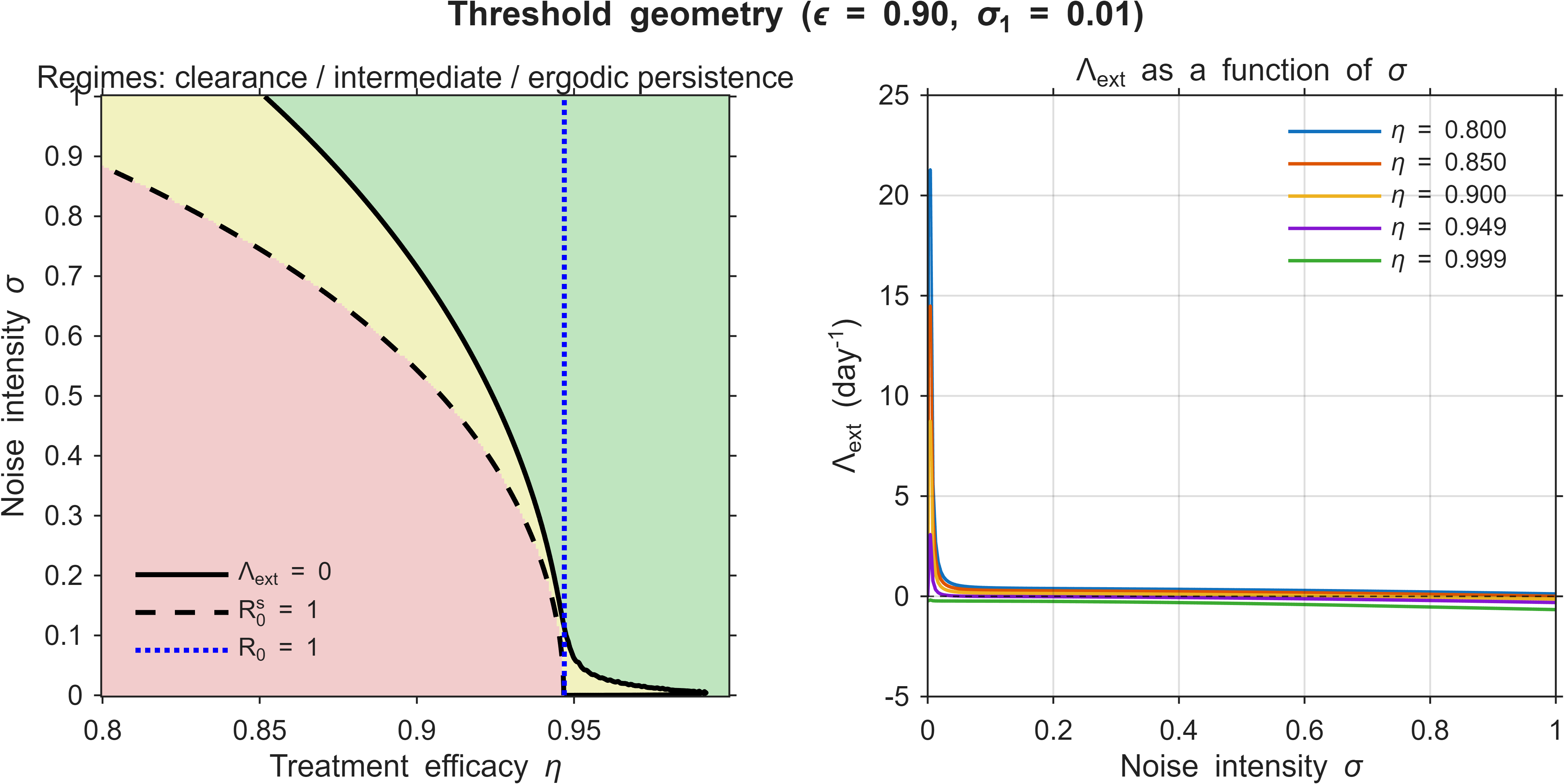}
\caption{\textbf{Extinction and persistence regions in the $(\eta,\sigma_3)$ plane}, with $\epsilon=0.90$, $\sigma_1=0.01$, and $\sigma_2=0.35$. Left: the reduced extinction boundary $\Lext=0$, the ergodic persistence boundary $\Rs=1$, and the deterministic threshold $\Rz=1$ are shown together with the intermediate region where neither extinction nor persistence is certified. Where condition $(\mathrm{C}_{\alpha^\ast})$ holds, $\Lext<0$ implies $\Lsharp=\Lext<0$ and therefore certifies almost sure exponential clearance by Theorem~\ref{Syz}. Right: $\Lext$ decreases with $\sigma_3$ for the selected values of $\eta$, illustrating the noise dependence established in Proposition~\ref{prop:detlimit}.}
\label{fig:regime}
\end{figure}

Figure~\ref{fig:regime} summarizes the theoretical regimes. Over the window shown, certified clearance occupies $39\%$ of the evaluated grid and certified persistence $52\%$, leaving an intermediate region where neither theorem is conclusive. Although the clearance region is plotted using $\Lext<0$, direct evaluation of $\Lsharp$ confirms that $\Lext<0$ implies $\Lsharp<0$ at every grid point in this region. The two bounds coincide except in a narrow strip near the boundary, where their difference is at most $3.3\times10^{-3}$ day$^{-1}$. Thus, every plotted clearance point is certified by Theorem~\ref{Syz}. The clearance boundary also shows that increasing $\sigma_3$ reduces the treatment efficacy required for clearance, but only gradually: decreasing $\eta$ from $0.94$ to $0.93$ increases the required virion noise from approximately $0.27$ to $0.44$. Thus, within the model, stochastic variability can partially offset reduced treatment efficacy, but at an unfavorable exchange rate.
\subsection{Long-Time Statistics, Discretization, and Phase-Space Geometry}

\begin{figure}[!htbp] 
\centering
\includegraphics[width=0.97\textwidth]{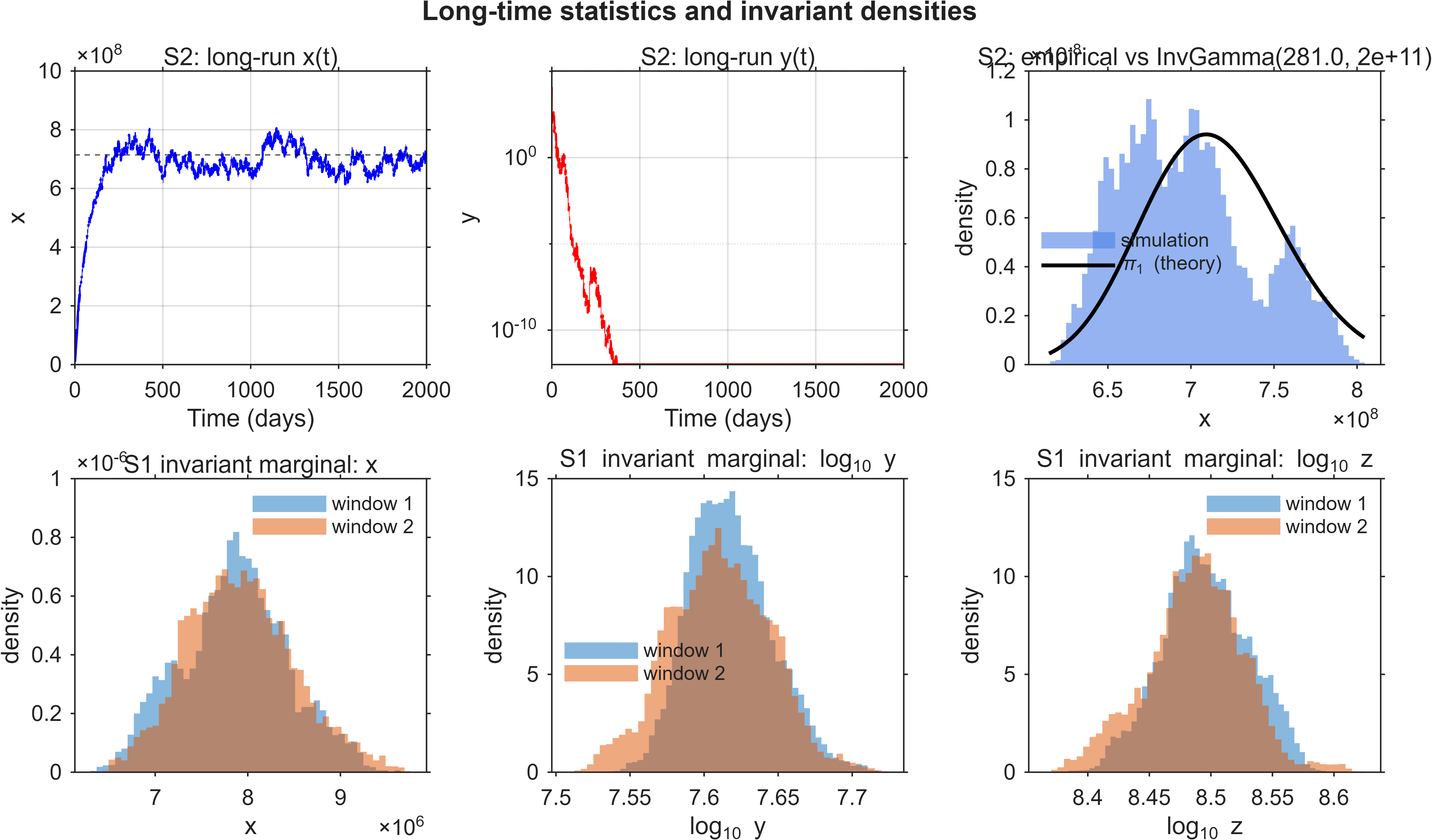}
\caption{\textbf{Long-time statistics and invariant distributions.} Top row: scenario S2 over $2000$ days, showing $x(t)$ relative to $\bar x$, $y(t)$ on a logarithmic scale, and the empirical density of $x$ compared with the inverse gamma law $\pi_1$ in \eqref{eq:invgamma}. Bottom row: scenario S1 over $2000$ days, showing empirical marginals of $x$, $\log_{10}y$, and $\log_{10}z$ over two disjoint terminal windows. Their close agreement is consistent with the ergodic invariant measure established in Theorem~\ref{thm:ergodic}.}
\label{fig:stationary}
\end{figure}

Figure~\ref{fig:stationary} tests the analysis rather than illustrating it. In the clearance regime, Theorem~\ref{x-x1} and Corollary~\ref{cor:attained} predict that the empirical law of $x$ converges to the explicit inverse gamma law $\pi_1$, whose second moment enters the extinction threshold through $\rho_1$ and, by Corollary~\ref{cor:attained}(ii), saturates the occupation bound \eqref{eq:secondmomentbound}. The script prints the sample mean and relative variance of $x$ for comparison with $\bar x$ and $\rho_1$; agreement to within a few percent is what a single trajectory of this length at $\Delta t=10^{-2}$ should give, the residual discrepancy reflecting time-discretization bias in the numerical invariant distribution together with finite-sample or Monte Carlo error. In the endemic regime, the marginals computed on disjoint terminal windows coincide, the signature one expects of an invariant measure. However, no finite simulation can establish that a stationary law exists or identify its support.

\begin{figure}[!htbp]
\centering
\includegraphics[width=0.52\textwidth]{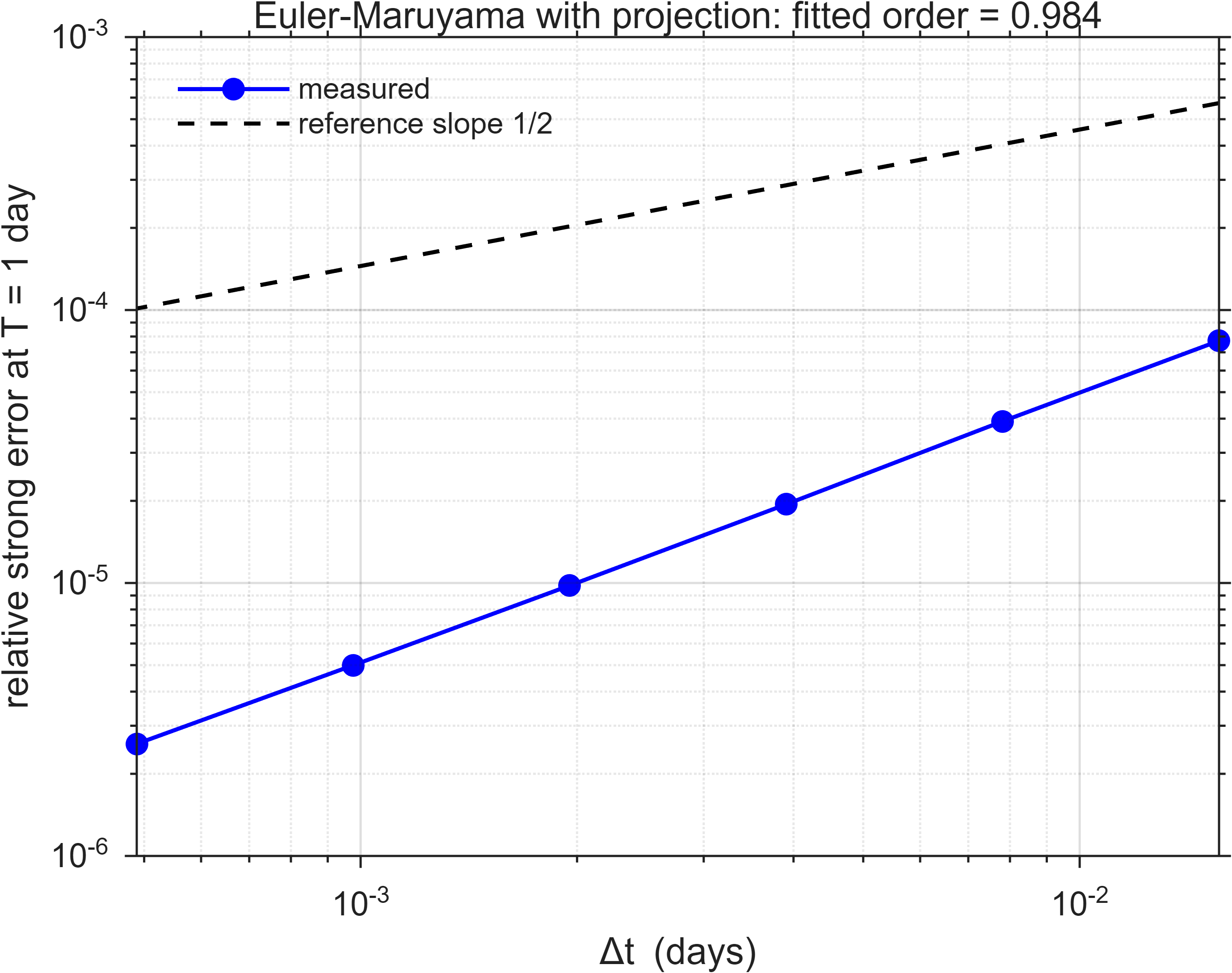}
\caption{\textbf{Strong-error convergence of the projected Euler--Maruyama scheme.} Relative strong error at $T=1$ day versus $\Delta t$, estimated from $200$ paths for scenario S1 using coupled Brownian increments and a reference step of $2^{-16}$ days. The dashed line has reference slope $1/2$; the measured log--log slope is $0.98$ over the tested range, indicating convergence at least consistent with the classical $1/2$ reference rate over this finite range.}
\label{fig:conv}
\end{figure}

Figure~\ref{fig:conv} quantifies the discretization error. Since no strong-order theorem is available for the projected scheme applied to \eqref{SM}, the experiment measures rather than verifies. The fitted slope, $0.98$, is not evidence of strong order one in the asymptotic sense: at the parameters of scenario S1 the noise intensities are small, $\sigma_i\le0.05$, while the drift is large, so over the range $\Delta t\in[2^{-11},2^{-6}]$ the error is dominated by the deterministic local truncation term, which is $O(\Delta t)$, and the stochastic component that carries the asymptotic $O(\Delta t^{1/2})$ behaviour has not yet emerged. The experiment supports the weaker, sufficient statement that the projection does not degrade the error below the classical reference rate over the tested range. The error at $\Delta t=10^{-3}$ lies several orders of magnitude below the features reported above, so the conclusions are not discretization artifacts. Two caveats remain. The projection biases the scheme upward near the boundary, so the far-tail decay rates of Figure~\ref{fig:noise} are rates for the discretization rather than for \eqref{SM}. And strong error at a fixed finite horizon says nothing about long-run statistics, which depend instead on how accurately the scheme reproduces the invariant distribution; we do not quantify that error for the projected scheme, and it is one possible source of the small discrepancy noted in Figure~\ref{fig:stationary}.

\begin{figure}[!htbp]
\centering
\includegraphics[width=0.85\textwidth]{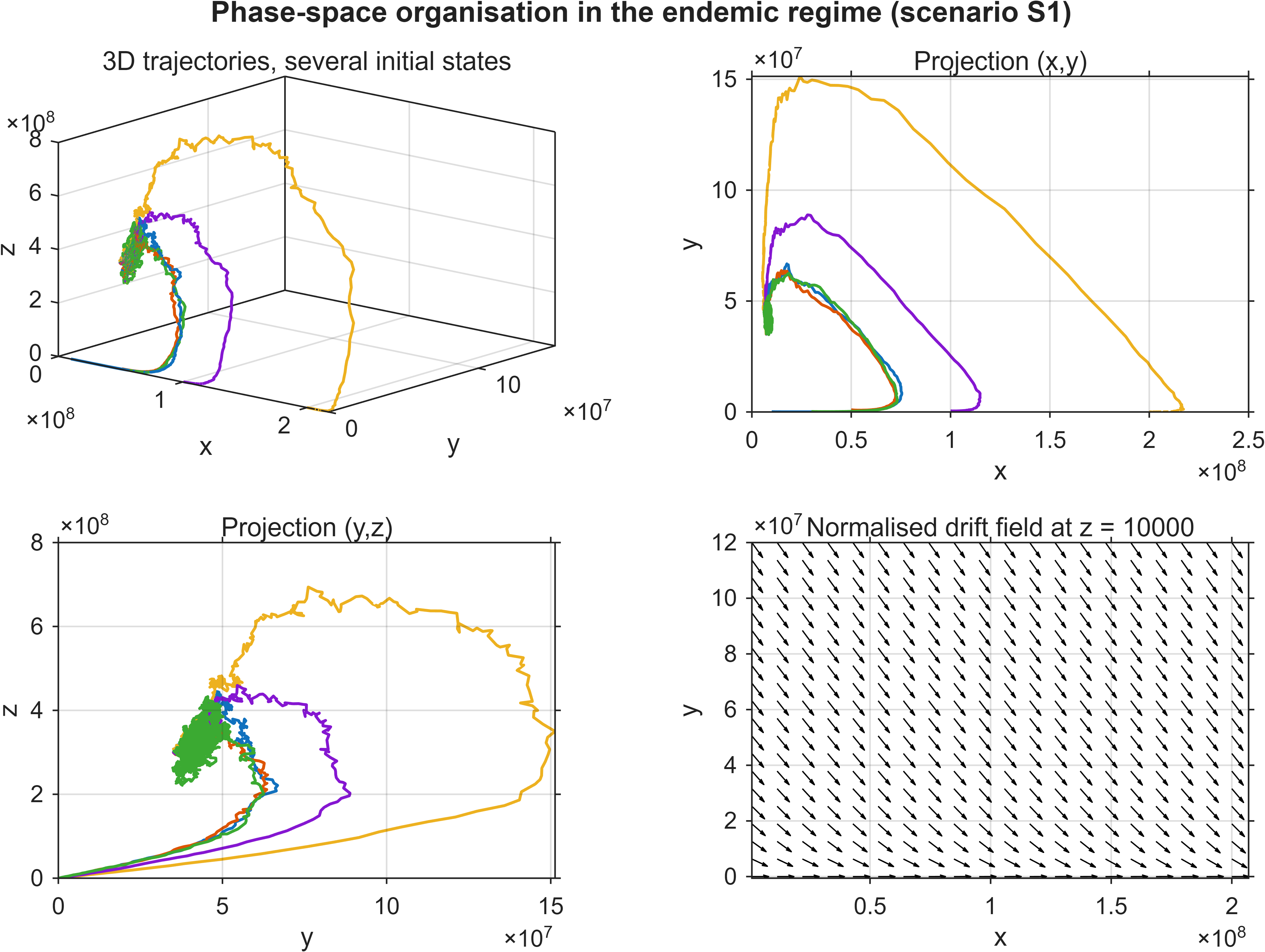}
\caption{\textbf{Phase-space organization in the endemic regime.} Top left: three-dimensional trajectories of scenario S1 launched from five different initial conditions under a common noise realization. Top right and bottom left: the $(x,y)$ and $(y,z)$ projections of the same trajectories. Bottom right: the normalized deterministic drift field of the $(x,y)$ subsystem at fixed $z=10^{4}$. Trajectories from distinct initial data leave their transients and concentrate in a common region of the interior, consistent with the existence of the invariant measure of Theorem~\ref{thm:ergodic}; the figure does not determine the support of that measure.}
\label{fig:phase}
\end{figure}

Figure~\ref{fig:phase} shows that a common long-run regime is reached from a wide range of initial conditions, the geometric counterpart of the drift condition \eqref{A2}: the Lyapunov function of Theorem~\ref{thm:ergodic} increases towards the boundary through the terms $-c_i\ln u_i$ and towards infinity through $W^{\theta+1}$, so excursions in either direction are penalized and the process returns to $D$. Recurrence of this kind does not imply that the invariant measure has compact support or is bounded away from the coordinate hyperplanes, and neither is claimed.

\section{Discussion}\label{sec:discussion}

The main mathematical contribution is a rigorous characterization of sufficient conditions for extinction and persistence in the stochastic HBV system. The weighted functional $\ln(\alpha y+z)$ improves substantially on the unweighted criterion and yields a reduced optimized bound that recovers the deterministic spectral threshold under the scaling of Definition~\ref{def:sharp}. A second key ingredient is the use of exact occupation identities for $x+y$. Because infection and cure cancel in this total cell population, these identities control hepatocyte fluctuations without replacing $x$ by its equilibrium or invoking an invalid pathwise comparison. This approach may be useful more broadly in compartmental stochastic models with internal transfer terms.

The biological interpretation should remain within the scope of the model. Stochastic fluctuations can shift the sufficient clearance condition, but the effect is primarily relevant near the deterministic threshold. In particular, parameter regimes with $\Rz>1$ may still satisfy the stochastic clearance criterion when the corresponding admissibility condition holds, whereas the noise levels required far from threshold are unrealistically large. Thus, the results describe a mathematical mechanism for noise-assisted clearance rather than suggesting that increased biological variability is therapeutically beneficial. The analysis also distinguishes clearance from persistence: under clearance, the hepatocyte population approaches a reference diffusion with an inverse gamma invariant law, whereas $\Rs>1$ yields a unique ergodic invariant distribution for the full system. Persistence is therefore represented by a stationary distribution rather than a deterministic endemic set point.

\subsection*{Limitations}

The stochastic forcing is multiplicative, white, and phenomenological, so the parameters $\sigma_i$ capture several sources of biological and measurement variability that the present model cannot distinguish. The framework also omits adaptive immunity, spatial structure, intracellular delays, and discrete infection events, and it does not attempt to infer patient-specific parameters. Extensions to correlated or jump-driven noise would provide more realistic representations of temporal and intermittent variability.

The theoretical criteria are sufficient, not necessary. The extinction result provides an upper bound on the almost sure Lyapunov exponent, with deterministic-limit sharpness established only under the scaling of Definition~\ref{def:sharp}, while the recurrence theorem leaves the intermediate region of Remark~\ref{rem:separation} unresolved. Hypothesis~(H2) supplies the fourth-moment control required for the martingale estimates, and condition $(\mathrm{C}_\alpha)$ arises from the occupation-measure relaxation used in the extinction proof. Removing these restrictions would require sharper pathwise or occupation-measure estimates. Finally, the projected numerical scheme preserves nonnegativity in the numerical sense. Still, it has no proved strong convergence order for this system and may reach the boundary, whereas the exact stochastic solution remains strictly positive at finite times.
\section{Conclusion}\label{sec:conclusion}

We developed a rigorous framework for the long-term dynamics of a stochastic within-host HBV model with multiplicative noise. The model admits a unique global positive solution with uniformly bounded moments. Using the weighted logarithmic functional $\ln(\alpha y+z)$, we derived an explicit sufficient condition for almost-sure exponential clearance and a reduced optimized bound $\Lext$. Under the admissibility condition $(\mathrm{C}_{\alpha^\ast})$, the general extinction bound $\Lsharp$ coincides with $\Lext$. Moreover, $\Lext$ satisfies
\[
\Lext \geq \lambda_{\max}(J)-\frac12\max\{\sigma_2^2,\sigma_3^2\},
\]
and converges to $\lambda_{\max}(J)$ under the deterministic-limit scaling $\sigma_2,\sigma_3\to0$ with $\sigma_1^2/(\sigma_2^2+\sigma_3^2)\to0$, thereby recovering the deterministic threshold $\Rz<1$. Exact occupation identities for the total hepatocyte population $x+y$ quantify the effect of target-cell fluctuations and yield a closed-form fluctuation penalty linked to the invariant inverse gamma law of the reference diffusion. Under clearance, the hepatocyte population converges almost surely and exponentially to this reference process. Conversely, when $\Rs>1$, the full stochastic system is positive recurrent and admits a unique ergodic invariant measure in the open positive orthant, these sufficient clearance and persistence regimes are mutually exclusive, leaving an intermediate parameter region unresolved. Numerical experiments support the analytical bounds, long-run distributions, and the conclusion that noise-induced clearance is primarily a near-threshold phenomenon.

Several extensions are natural. Correlated or jump-driven noise could provide more realistic representations of immune and treatment variability. Resolving the intermediate regime may require quasi-stationary or boundary-asymptotic methods. Finally, treating $\eta$ and $\epsilon$ as control variables would enable treatment optimization subject to the rigorous stochastic clearance condition $\Lsharp<0$, rather than the deterministic criterion $\Rz<1$.


\section*{Conflict of Interest}

The author declares no conflict of interest.

\end{document}